\documentclass[12pt]{article}

\usepackage[T1]{fontenc}
\usepackage[latin2]{inputenc}
\usepackage{lmodern}
\usepackage{amssymb}
\usepackage{stmaryrd}
\usepackage{theorem}

\usepackage{amsmath}
\usepackage{amsfonts}
\usepackage{graphicx}
\usepackage{epstopdf}
\usepackage[parfill]{parskip}    
\usepackage{hyperref}
\hypersetup{pdftex,colorlinks=true,allcolors=blue}
\usepackage{hypcap}
\usepackage{algorithm}
\usepackage{algpseudocode} 
\usepackage{authblk}
\usepackage{tabto}

\newtheorem{thm}{Theorem}[section]
\newtheorem{dfn}{Definition}[section]

\newtheorem{prm}{Problem}[section]

\newcommand*{\sm}{\ensuremath{\sim}}
\newcommand*{\st}{\ensuremath{\otimes}}
\newcommand*{\co}{\ensuremath{\circ}}
\newcommand*{\sd}{\ensuremath{\triangle}}
\newcommand*{\be}{\ensuremath{\bar{e}}}
\newcommand*{\beo}{\ensuremath{\be_1}}
\newcommand*{\bet}{\ensuremath{\be_2}}
\newcommand*{\beh}{\ensuremath{\be_3}}
\newcommand*{\im}{\ensuremath{\in_m}}
\newcommand*{\imo}{\ensuremath{\in_{m_1}}}
\newcommand*{\imt}{\ensuremath{\in_{m_2}}}
\newcommand*{\N}{\ensuremath{\mathbb{N}}}
\newcommand*{\Z}{\ensuremath{\mathbb{Z}}}
\newcommand*{\D}{\ensuremath{\mathrm{D}}}
\newcommand*{\R}{\ensuremath{\mathrm{R}}}
\newcommand*{\E}{\ensuremath{\mathrm{E}}}
\newcommand*{\M}{\ensuremath{\mathrm{M}}}
\newcommand*{\B}{\ensuremath{\mathrm{B}}}
\newcommand*{\Bm}{\ensuremath{\mathrm{B}_m}}
\newcommand*{\Sy}{\ensuremath{\mathrm{S}}}
\newcommand*{\Sm}{\ensuremath{\mathrm{S}_m}}
\newcommand*{\Zm}{\ensuremath{\mathbb{Z}_m}}
\newcommand*{\Nm}{\ensuremath{\mathrm{N}_m}}
\newcommand*{\Nmo}{\ensuremath{\mathrm{N}_{m_1}}}
\newcommand*{\Nmt}{\ensuremath{\mathrm{N}_{m_2}}}
\newcommand*{\Nme}{\ensuremath{\mathrm{N}_m^e}}
\newcommand*{\Rm}{\ensuremath{\mathrm{R}_m}}
\newcommand*{\kRm}{\ensuremath{_k\mathrm{R}_m}}
\newcommand*{\kRme}{\ensuremath{\kRm^e}}
\newcommand*{\Rmo}{\ensuremath{\mathrm{R}_{m_1}}}
\newcommand*{\Rmt}{\ensuremath{\mathrm{R}_{m_2}}}
\newcommand*{\Rme}{\ensuremath{\mathrm{R}_m^e}}
\newcommand*{\Gm}{\ensuremath{\mathrm{G}_m}}
\newcommand*{\Gmo}{\ensuremath{\mathrm{G}_{m_1}}}
\newcommand*{\Gmt}{\ensuremath{\mathrm{G}_{m_2}}}

\newcommand*{\Em}{\ensuremath{\mathrm{E}_m}}
\newcommand*{\Emo}{\ensuremath{\mathrm{E}_{m_1}}}
\newcommand*{\Emt}{\ensuremath{\mathrm{E}_{m_2}}}
\newcommand*{\prf}{\textbf{Proof}\ \ }
\newcommand*{\sln}{\textbf{Solution}\ \ }
\newcommand*{\sqr}{\ensuremath{\square}}
\newcommand*{\rar}{\ensuremath{\Rightarrow}}
\newcommand*{\lar}{\ensuremath{\Leftarrow}}
\newcommand*{\lrar}{\ensuremath{\Leftrightarrow}}
\newcommand*{\phm}{\ensuremath{\varphi(m)}}
\newcommand*{\psm}{\ensuremath{\psi(m)}}
\newcommand*{\am}{\ensuremath{|a|_m}}
\newcommand*{\bm}{\ensuremath{|b|_m}}
\newcommand*{\cm}{\ensuremath{|c|_m}}
\newcommand*{\ibma}{\ensuremath{\mathrm{ind}_b^m a}}

\newcommand*{\indm}[2]{\ensuremath{\mathrm{ind}_{#1}^m #2}}
\newcommand*{\eibma}{\ensuremath{\exists\mathrm{ind}_b^m a}}

\newcommand*{\orbm}[1]{\ensuremath{\langle #1 \rangle_m}}
\newcommand*{\oma}{\ensuremath{\omega_m (a)}}
\newcommand*{\Oma}{\ensuremath{\Omega_m (a)}}
\newcommand*{\om}{\ensuremath{\omega(m)}}
\newcommand*{\ome}[1]{\ensuremath{\omega\left(#1\right)}}
\newcommand*{\ind}{\ensuremath{\mathrm{ind}}}
\newcommand*{\ali}{\ensuremath{\alpha_i}}
\newcommand*{\piai}{\ensuremath{p_i^{\alpha_i}}}
\newcommand*{\pibi}{\ensuremath{p_i^{\beta_i}}}
\newcommand*{\mmod}{\ensuremath{\mathrm{mod\ }}}
\newcommand*{\Db}{\ensuremath{\D_m(b,c)}}
\newcommand*{\Dc}{\ensuremath{\D_m(c,b)}}
\newcommand*{\mand}{\ensuremath{\mathrm{and}}}
\newcommand*{\mor}{\ensuremath{\mathrm{or}}}
\newcommand*{\Mmka}{\ensuremath{\mathrm{M}_m(k,a)}}

\newcommand*{\SmRka}{\ensuremath{\mathrm{S}_m^{\mathrm{R}}(k,a)}}
\newcommand*{\SmR}{\ensuremath{\mathrm{S}_m^{\mathrm{R}}}}
\newcommand*{\xni}{\ensuremath{x_0^{-1}}}
\newcommand*{\mma}{\ensuremath{\mu_m(a)}}

\newcommand*{\rhme}{\ensuremath{\rho_m^e}}
\newcommand*{\rhmo}{\ensuremath{\rho_m^1}}
\newcommand*{\rme}{\ensuremath{r_m^e}}
\newcommand*{\rmo}{\ensuremath{r_m^1}}
\newcommand*{\eqm}{\ensuremath{\sim_m}}
\newcommand*{\QM}{\ensuremath{\mathcal{QM}}}
\newcommand*{\DI}{\ensuremath{\mathcal{DI}}}
\newcommand*{\DIpa}{\ensuremath{\mathcal{DI}_{p^{\alpha}}}}
\newcommand*{\MU}{\ensuremath{\mathcal{M}}}
\newcommand*{\Domf}{\ensuremath{\mathcal{D}(f)}}
\newcommand*{\Dom}{\ensuremath{\mathcal{D}}}
\newcommand*{\lcm}{\ensuremath{\mathrm{lcm}}}
\newcommand*{\Smk}{\ensuremath{\mathrm{S}_{m,k}}}
\newcommand*{\br}{\ensuremath{\bar{r}}}

\begin{document}


\begin{titlepage}

\begin{center}

\vspace*{2cm}

\Large
On Composite Moduli from the Viewpoint of\\
Idempotent Numbers

\vspace{1cm}
\large
MSc Thesis

\vspace{1cm}
József Vass\\
\small
\vspace{0.25cm}
jozsef.vass@outlook.com


\begin{figure}[H]
\centering
\vspace{7cm}
\includegraphics[width=85pt]{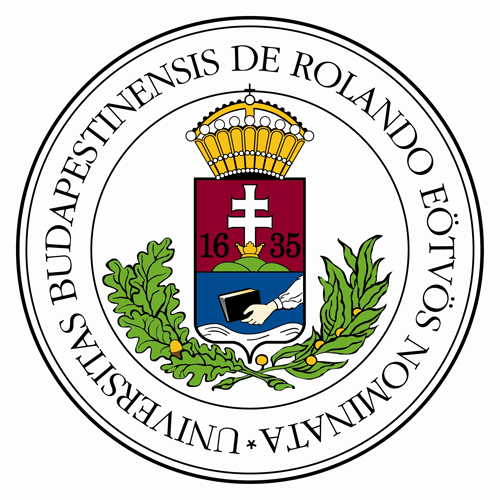}
\end{figure}

E\"otv\"os Lor\'and University\\
Budapest, Hungary\\
2004

\normalsize

\end{center}

\end{titlepage}

\tableofcontents

\newpage

\section{Introduction}

The purpose of this paper is to introduce basic concepts, as of yet unknown, that are fundamental in the examination of composite moduli, while avoiding the notoriously difficult problem of prime-factorization. We introduce a new class of numbers, called the set of idempotent numbers, that is unavoidable when researching composite moduli. Among many interesting results, we give generalizations of well-known theorems and definitions, such as the Euler-Fermat Theorem and the concept of primitive roots. We consider the generalization of the equivalence condition for the solvability of a binomial congruence to be the main result of our paper.

The paper is organized into sections of definitions, theorems, and notes. We intended to include every single result known to us regarding idempotent numbers, so as to propagate any further research that may be done on the subject by other authors. We provide basically no references, since we were unable to find any, on this particular part of Number Theory. Some of the simpler results however, may be known to group theorists, or may be found among exercises in textbooks. Since we publish all of our results, it is no surprise that the reader may find incomplete lines of thought within the network of theorems. We welcome and surely appreciate any helpful comments and papers that would be related to them.

Many of our notes in this paper are meant to provide insight to the reader into our aims of research currently in progress. As we have formerly mentioned, the paper itself is concerned with revealing all the known implications of the existence of a set, called the set of idempotent numbers, modulo a composite number.

In the next section, we provide definitions and notations, which will be used throughout the paper.

In the section ``Idempotent Numbers'', we give results which show the fundamentality of this set, whenever we wish to explore
the hazy structure of composite moduli. We also give a generalization of the definition of order, based on our generalization of the well-known Euler-Fermat Theorem.

In the section ``Normal Numbers'', we define a set, which has quite a hazy structure, although not as much as the whole of $\Zm$.

In the section ``Regular Numbers'', we introduce the nicest and most general set, that one can work with, when examining composite moduli. We are going to show many properties for this set, that have only been known so far for its subset of reduced residues, and we will also partition it into subsets of Abelian groups. We also give several other definitions and results, to the best of our knowledge, as of yet unknown.

The next section is the most important one, since it contains the main result of this paper, which is an equivalence condition for the solvability of a binomial congruence, to the greatest degree of generality, we could hope to reach. This condition is currently difficult to calculate in practice.

The next section generalizes the well-known definition of primitive roots to composite moduli. It includes only a few results, which implies that it needs the most intensive level of research, since we believe that revealing the structure of generalized primitive roots, shall also solve the problem of calculating the condition we have spoken of above.

The section ``Number Theoretic Functions'', introduces functions defined through definitions of previous sections, all of them
containing the idea of idempotent numbers.

The set of idempotent numbers shows interesting algebraic properties, as revealed in the next section. Some nice operators may
be defined over its elements, which are analogous to operators known from Set Theory. We will also show that
idempotent numbers form a commutative ring.

Our last section, entitled ``Second-Degree Polynomials'', discusses their sets of solutions, and in some special cases, characterizes them as well, with the use of idempotent numbers.

\newpage

\section{Basic Definitions}

Let \N\ denote the set of whole numbers greater than or equal to $1$. By ``number'' we will mean any whole number. In the entire paper, let $m$ denote a fixed integer, and let
\[ m=\prod_{i=1}^{\infty} \piai,\ p_i\ \mathrm{prime},\ p_i<p_{i+1},\ \alpha_i\in\N\cup\{0\},\ i\in\N \]
be the prime-factorization of $m$. This $m$ is said to be square-free if $\alpha_i\le 1$ for all $i\in\N$.
\[ \om := |\{i\in\N: \alpha_i>0\}|. \]
We shall call $m$ a weakly even number, if $p_1=2$ and $0\le\alpha_1\le 2$. Furthermore, we shall call $m$ a barely even number, if $p_1=2$ and $\alpha_1=1$.

Let $\rar$ denote logical inference. Let ``iff'' mean logical equivalence, or ``if and only if'', and let it be denoted by $\lrar$. Let $\exists$ denote existence ($\exists!$ meaning ``exists exactly one''), and $\forall$ the logical term ``for all''.

Denote $\Zm := \{1,\dots,m\}$. In case of $A\subset\Zm,\ a\in\Z$
\[ a\im A \lrar \exists b\in A: a\equiv b\ (\mmod m). \]
Let $a\ \mmod m$ denote the number $b\in\Zm$ for which $a\equiv b\ (\mmod m)$. In case of $A\subset\Z$
\[ A\ \mmod m := \{a\ \mmod m:\ a\in A\}. \]

Let $(a,b)$ denote the greatest common divisor of the numbers $a$ and $b$. Also, in case of $A\subset\N$, let $\mathrm{gcd} (a: a\in A)$ denote the greatest common divisor of all the elements in $A$. Let $[a,b]$ denote the least common multiple of the numbers $a$ and $b$. Also, in case of $A\subset\N$, let $\mathrm{lcm} (a: a\in A)$ denote the least common multiple of all the
elements in $A$. In case of an integral vector $a\in\N^n,\ n\in\N$, let $[a]$ denote the least common multiple of the coordinates of $a$. Let \phm\ denote the number of integers relatively prime to, and not exceeding $m$. Furthermore, let \psm\ denote
\[ \psm:= \mathrm{lcm} (\varphi(\piai):\ \alpha_i>0). \]
(Note that if $m$ is weakly even, then the maximal order modulo $m$ is \psm.)

For $a,b\in\Z$, let $a\mid b$ denote that $a$ is a divisor of $b$. In case of vectors $a,b\in\Z^n,\ n\in\N$, and any relation $\sim$, let $a\sim b$ mean that $a_i\sim b_i$ for all $1\le i\le n$. Let $Dom(f)$ denote the domain of function $f$.

Regarding other basic notations and theorems on congruences, see \cite{bb00020}.

\newpage
\section{Idempotent Numbers}

\begin{dfn} \label{in01}
A number $e$ is said to be idempotent modulo m, if
\[ e^2\equiv e\ (\mmod m). \]
Let \Em\ denote the subset of idempotent numbers in \Zm.
\end{dfn}

We will mostly denote an idempotent number by $e$. The notation comes from the first letter of the Hungarian word for ``unit'', since there may be defined groups in \Zm, with their units being idempotent numbers modulo $m$. If $m$ has only one prime factor, then $\Em=\{1,m\}$.

\begin{thm} \label{in02}
For all $a\in\Z$
\[ a^{\phm} \im \Em. \]
\end{thm}
\noindent
\prf
Let $i\in\N$ be a fixed number, and $\alpha_i>0$. There are two possibilities.\\
1. $p_i\mid a$\\
Since
\[ \alpha_i = 1+(\alpha_i -1)\le 2^{\alpha_i -1}\le p_i^{\alpha_i -1}\le p_i^{\alpha_i -1}(p_i -1)
\le \phm \]
we have
\[ a^{\phm}\equiv 0\ (\mmod \piai). \]
2. $p_i\nmid a$\\
In this case, by the Euler-Fermat Theorem (and $\varphi (p_i^{\alpha_i})\mid\phm$), we have
\[ a^{\phm}\equiv 1\ (\mmod \piai). \]
In both cases 1. and 2., we have
\[ a^{\phm}(a^{\phm}-1)\equiv 0\ (\mmod \piai) \]
so for all $i\in\N$ we have
\[ (a^{\phm})^2\equiv a^{\phm}\ (\mmod \piai) \]
which means that
\[ (a^{\phm})^2\equiv a^{\phm}\ (\mmod m) \]
so $a^{\phm}$ is idempotent modulo $m$. \sqr

\begin{thm} \label{in12}
For all $a\in\Z$
\[ a^{\psm} \im \Em. \]
\end{thm}
\noindent
\prf
The proof goes the same way as above, by changing each \phm\ to \psm. \sqr

Many interesting facts follow from the two theorems above. One is that every polynomial is equivalent to a
polynomial of degree less than or equal to $2\psm -1$ (modulo $m$).

\begin{thm} \label{in08}
For all $a\in\Z$
\[ a^{\phm}\equiv (a,m)^{\phm}\ (\mmod m). \]
\end{thm}
\noindent
\prf
This fact may be proven the same way as Theorem \ref{in02}, by considering that
\[ p_i\mid a\ \lrar\ p_i\mid (a,m)\ \ (\alpha_i>0).\ \ \sqr \]

\begin{thm} \label{in06}
\[ |\Em|=2^{\om}. \]
\end{thm}
\noindent
\prf
From the proof of Theorem \ref{in02}, we see that a number is idempotent modulo $m$ iff it
is congruent to either $0$ or $1$ modulo each of the prime power divisors of $m$.
From this fact, our theorem follows quite clearly, by the application of the Chinese
Remainder Theorem (see \cite{bb00020}). \sqr

\begin{thm} \label{in07}
For $k\in\N$
\[ |k\Em\ \mmod m|=2^{\ome{\frac{m}{(k,m)}}}. \]
\end{thm}
\noindent
\prf
Let
\[ \frac{m}{(k,m)}=\prod_{i=1}^{\infty} \pibi \]
be the prime-factorization of $\frac{m}{(k,m)}$. We see that $\beta_i=0$ iff $\piai\mid k$, so
\[ \ome{\frac{m}{(k,m)}} = |\{i:\beta_i>0\}|=|\{i:\piai\nmid k\}|. \]
Without hurting generality, we may suppose that $p_i$ are ordered so that for some $n\in\N$ we have $\beta_i>0$ if
$1\le i\le n$ and $\beta_i=0$ if $i>n$. Let
\[ m_1:=\prod_{i=1}^n \piai,\ \ m_2:= \prod_{i=n+1}^{\infty} \piai. \]
This way we have for all $e\in\Em$ that
\[ ke\equiv 0\ \mor\ ke\equiv k \not\equiv 0\ (\mmod \piai)\ \ (1\le i\le n) \]
and $ke\equiv 0\ (\mmod m_2)$. So by the Chinese Remainder Theorem and Theorem \ref{in06} we have
\[ |k\Em\ \mmod m|=|k\Emo\ \mmod m_1|\cdot |k\Emt\ \mmod m_2|=|\Emo|\cdot |\{ 0\}|=2^n=2^{\ome{\frac{m}{(k,m)}}}.\ \ \sqr \]

\begin{thm} \label{in03}
If for some $a\in\Z,\ k,l\in\N$, we have $a^k, a^l \im\Em$, then $a^k\equiv a^l\ (\mmod m)$.
\end{thm}
\noindent
\prf
\[ a^k\equiv (a^k)^l\equiv a^{kl}\equiv (a^l)^k\equiv a^l\ (\mmod m).\ \ \sqr \]

\begin{dfn} \label{in04}
For $a\in\Z$, let its order modulo m to be the smallest $n\in\N$ for which $a^n \im \Em$. Let \am\ denote this $n$.
\[ a^0:=a^{\am}\ \mmod m. \]
Furthermore, let the inverse of $a$ be denoted as
\[ a^{-1}:= a^{\am -1}\ \mmod m \]
and for $k\in\N$
\[ a^{-k}:= (a^{-1})^k\ \mmod m. \]
For $b\in\Z$, if it exists, let \ibma\ denote the smallest $k\in\N$, for which $b^k\equiv a\ (\mmod
m)$, and let its existence be denoted as \eibma.
\end{dfn}

The existence of the above $n$ is guaranteed by Theorem \ref{in02}.

\begin{thm} \label{in11}
If for some $a\in\Z,\ k,n\in\N,\ k\le n$ we have $a^{k+n}\equiv a^k\ (\mmod m)$, then $a^n\im\Em$.
\end{thm}
\noindent
\prf
\[ a^{k+n}\equiv a^k\ \rar\ a^{n-k}a^{k+n}\equiv a^{n-k}a^k\ \rar\ (a^n)^2\equiv a^n\ (\mmod m) \]
so we have that $a^n\im\Em$. \sqr

\begin{thm} \label{in05}
If $m_1,m_2\in\N,\ m=[m_1,m_2]$, then
\[ e\in\Em\ \lrar\ e\imo\Emo\ \mathrm{and}\ e\imt\Emt. \]
\end{thm}
\noindent
\prf
This follows from the fact that a number is idempotent modulo $m$ iff it is congruent to either $0$
or $1$ modulo each of the prime power divisors of $m$. \sqr

\begin{prm} \label{in10}
Defining the sequence of numbers $(a_n)_{n=1}^{\infty}$ recursively with
\[ a_0:=1,\ a_{n}:=42^{a_{n-1}}\ (n\in\N) \]
what will be the last two digits of $a_{100}$?
\end{prm}
\noindent
\sln
With some calculation, we get the following results
\[ a_{100}=42^{a_{99}},\ |42|_{100}=20,\ 42^{20}\equiv 76\ (\mmod 100) \]
\[ a_{99}=42^{a_{98}},\ |42|_{20}=4,\ 42^4\equiv 16\ (\mmod 20) \]
\[ a_{98}=42^{a_{97}},\ |42|_4=2,\ 42^2\equiv 4\ (\mmod 4) \]
\[ a_{97}\equiv 0\ (\mmod 2)\ \rar\ a_{98}\equiv 0\ (\mmod 4)\ \rar\ a_{99}\equiv 16\ (\mmod 20)\ \rar \]
\[ \rar a_{100}\equiv 76\cdot 42^{16}\equiv 76\cdot 56\equiv 56\ (\mmod 100).\ \ \sqr \]

Note that the above problem may be solved in other ways as well. I decided to include it in our discussion, because this idea of a solution made me realize the importance of idempotent numbers when discussing composite moduli, and it also gave me incentive to investigate composite moduli from the viewpoint of idempotent numbers, starting the research which resulted in this paper.

\newpage
\section{Normal Numbers}

\begin{dfn} \label{nn01}
$a\in\Z$ is said to be normal modulo $m$ if the following logical inference holds
\[ a^k\im\Em\ \rar\ \am\mid k\ \ (k\in\N). \]
Let \Nm\ denote the subset of normal numbers in \Zm. Furthermore, for $e\in\Em$, let \Nme\ denote the set
\[ \{a\in\Nm:\ a^{\am}\equiv e\ (\mmod m)\}. \]
\end{dfn}

\begin{thm} \label{nn02}
A number $a\in\Zm$ is normal iff the following logical inference holds
\[ a^k\equiv a^l\ (\mmod m)\ \rar\ k\equiv l\ (\mmod \am)\ \ (k,l\in\N). \]
\end{thm}
\noindent
\prf
First, let us suppose that $a\in\Zm$ is normal. Then
\[ a^k\equiv a^l\ \rar\ a^k a^{l\phm -l}\equiv a^{l\phm}\ (\mmod m). \]
Since $a^{l\phm}\im\Em$ and $\am\mid\phm$, we have
\[ 0\equiv k+l\phm -l\equiv k-l\ \rar\ k\equiv l\ (\mmod \am). \]
Now, if the inference holds, with $l:=\am$ we have that $a$ is normal. \sqr

\begin{thm} \label{nn06}
If $m_1,m_2\in\N,\ m_1\mid m_2,\ a\in\Nmo$, then $|a|_{m_1}\mid |a|_{m_2}$.
\end{thm}
\noindent
\prf
\[ a^{|a|_{m_2}}\imt\Emt,\ m_1\mid m_2\ \rar\ a^{|a|_{m_2}}\imo\Emo,\ a\in\Nmo\ \rar\ |a|_{m_1}\mid |a|_{m_2}.\ \ \sqr \]

\begin{thm} \label{nn04}
Let $m_1,m_2\in\N$ be such that $m=[m_1,m_2]$, and $a\in\Zm$. If $a\imo\Nmo$ and $a\imt\Nmt$, then $a\in\Nm$ and $\am=[|a|_{m_1},|a|_{m_2}]$.
\end{thm}
\noindent
\prf
By Theorem \ref{in05}, for $k\in\N$ we have
\[ a^k\im\Em\ \lrar a^k\imo\Emo,\ a^k\imt\Emt\ \rar\ [|a|_{m_1},|a|_{m_2}]\mid k\ \]
so $[|a|_{m_1},|a|_{m_2}]\mid \am$. With $k=[|a|_{m_1},|a|_{m_2}]$ we have $a^{[|a|_{m_1},|a|_{m_2}]}\im\Em$ so
$\am\le [|a|_{m_1},|a|_{m_2}]$, so $\am=[|a|_{m_1},|a|_{m_2}]$. We also see from above that $a$ is normal modulo $m$. \sqr

\begin{thm} \label{nn03}
If $a\in\Nm,\ k\in\N$ then
\[ |a^k|_m = \frac{\am}{(k,\am)}. \]
\end{thm}
\noindent
\prf
\[ (a^k)^{\frac{\am}{(k,\am)}}\equiv (a^{\am})^{\frac{k}{(k,\am)}}\im\Em \]
so we have
\[ |a^k|_m\le \frac{\am}{(k,\am)} \]
by Definition \ref{in04}.
\[ a^{kl}\im\Em\ \rar\ \am\mid kl\ \lrar\ \frac{\am}{(k,\am)}\mid l \]
so we have
\[ |a^k|_m\ge \frac{\am}{(k,\am)}.\ \ \sqr \]

\begin{thm} \label{nn05}
If $a\in\Nm,\ n\in\N$, then $a^n\im\Nm$.
\end{thm}
\noindent
\prf
For all $k\in\N$
\[ (a^n)^k\im\Em\ \rar\ \am\mid nk\ \lrar\ |a^n|_m=\frac{\am}{(n,\am)}\mid k.\ \ \sqr \]

\begin{thm} \label{nn07}
Let $e\in\Em,\ a,b\in\Nme,\ k\in\N$ be such that $\eibma$. Then
\[ a^{\frac{\bm}{(k,\bm)}}\im\Em\ \rar\ (k,\bm)\mid\ibma. \]
\end{thm}
\noindent
\prf
\[ e\equiv a^{\frac{\bm}{(k,\bm)}}\equiv b^{\frac{\bm\ibma}{(k,\bm)}}\ (\mmod m)\ \rar \]
\[ \rar\ \bm\mid \frac{\bm\ibma}{(k,\bm)}\ \lrar\ \frac{\ibma}{(k,\bm)}\in\Z\ \lrar\ (k,\bm)\mid\ibma.\ \ \sqr \]

\begin{thm} \label{nn08}
For all $a\in\Nm$ we have $a^{-1}\in\Nm,\ |a^{-1}|_m = |a|_m$ and
\[ (a^{-1})^{-1}\equiv a^{\am +1}\ (\mmod m). \]
\end{thm}
\noindent
\prf
The first statement follows from Theorem \ref{nn05}. The second statement follows from Theorems \ref{nn05} and \ref{nn03} since
\[ |a^{\am -1}|_m = \frac{\am}{(\am-1,\am)} = \am. \]
If $\am\ge 2$ then
\[ (a^{-1})^{-1}\equiv (a^{\am -1})^{\am -1}\equiv a\cdot (a^{\am})^{\am -2}\ (\mmod m). \]
The case of $a\in\Em$ is trivial. \sqr

\newpage
\section{Regular Numbers}

\begin{dfn} \label{rn01}
$a\in\Z$ is said to be regular modulo $m$ if
\[ a^{\am +1}\equiv a\ (\mmod m). \]
Let \Rm\ denote the subset of regular numbers in \Zm. Furthermore, for $e\in\Em$, let \Rme\ denote the set
\[ \{a\in\Rm:\ a^{\am}\equiv e\ (\mmod m)\}. \]
\end{dfn}

From Theorems \ref{in02} and \ref{in03} we have that for all $a\in\Z$
\[ a^{\am}\equiv a^{\phm}\ (\mmod m). \]
We will apply this simple fact in our next theorem.

\begin{thm} \label{rn19}
$\Rm=\Zm$ iff $m$ is square-free.
\end{thm}
\noindent
\prf
First, let us suppose that $\Rm=\Zm$. Then for all $i:\alpha_i>0$ we have
\[ p_i(p_i^{\phm}-1)\equiv 0\ (\mmod m)\ \rar\ p_i(p_i^{\phm} -1)\equiv 0\ (\mmod \piai)\ \rar\ \alpha_i=1. \]
Now, if $m$ is square-free, then for all $a\in\Zm$
\[ a(a^{\phm}-1)\equiv 0\ (\mmod p_i)\ (i\in\N)\ \rar\ a(a^{\phm}-1)\equiv 0\ (\mmod m).\ \ \sqr \]

\begin{thm} \label{rn02}
\[ \Rm \subset \Nm. \]
\end{thm}
\noindent
\prf
We need to show that for all $a\in\Rm$, if $a^k\im\Em$ then $\am\mid k$ for all $k\in\N$. Let $q,r\in\N\cup\{0\}$ be such that
\[ k=q\am +r,\ 0\le r<\am. \]
Let us suppose that $r>0$. Then we have
\[ a^r\equiv a\cdot a^{r-1}\equiv a^{\am +1}\cdot a^{r-1}\equiv (a^{\am})^q\cdot a^r\equiv a^k\ (\mmod m) \]
so $a^r\im\Em$, which is a contradiction by Definition \ref{in04}. \sqr

\begin{thm} \label{rn31}
For all $e\in\Em$
\[ \Rme = \{ea\ \mmod m:\ a^{\am}\equiv e\ (\mmod m),\ a\in\Zm \}. \]
\end{thm}
\noindent
\prf
If $a\in\Rme$ then $a\equiv e\cdot a\ (\mmod m)$ which is obviously in the set on the right-hand side. Now if $a\in\Zm$ is such that $a^{\phm}\equiv e\ (\mmod m)$, then multiplying this congruence by $ea$ we have
\[ (ea)(ea)^{\phm}\equiv ea\ (\mmod m) \]
so $ea\im\Rme$. \sqr

\begin{thm} \label{rn18}
For all $a\in\Z$
\[ a^m\equiv a^{m-\phm}\equiv a^{m+\phm}\ (\mmod m) \]
and
\[ a^{m-\phm}\im\Rm. \]
Furthermore
\[ a^{\phm -1}\im\Rm \]
if $\am<\phm$ and $a\in\Nm$.
\end{thm}
\noindent
\prf
Let $i\in\N$ be a fixed number, and $\alpha_i>0$. There are two possibilities.\\
1. $p_i\mid a$
\[ \alpha_i\le p_i^{\alpha_i -1}\le p_i^{\alpha_i -1}\cdot q = m-\phm \]
for some $q\in\N$, so
\[ a^{m-\phm}\equiv 0\ (\mmod \piai). \]
2. $p_i\nmid a$\\
In this case we have
\[ a^{\phm}\equiv 1\ (\mmod \piai). \]
In both cases 1. and 2., we have
\[ a^{m-\phm}(a^{\phm} -1)\equiv 0\ (\mmod \piai) \]
from which we have our first congruence. The second one follows easily, by multiplying both sides by $a^{\phm}$.
\[ a^{m-\phm}\cdot (a^{m-\phm})^{\phm}\equiv a^{m-\phm}\cdot a^{\phm}=a^m\equiv a^{m-\phm}\ (\mmod m). \]
From Theorem \ref{rn02}, and the definition of normal numbers, we have that $\phm=k\am$ for some $k>1$. So
\[ a^{\phm -1}\cdot (a^{\phm -1})^{\phm}\equiv a^{\phm -1}\cdot a^{\phm} \equiv \]
\[ \equiv a^{\phm}\cdot a^{\am}\cdot a^{(k-1)\am -1} \equiv a^{\am}\cdot a^{(k-1)\am -1} = a^{\phm-1}\ (\mmod m).\ \ \sqr \]

Note that based on our previous theorem, we may define the function
\[ \delta_m(a) := \min \{n\in\N:\ a^n\im\Rm\}\ \ (a\in\Zm) \]
which has the property
\[ \delta_m(a) = \min \{n\in\N:\ a^{\am +n}\equiv a^n\ (\mmod m)\} \]
since
\[ a^n\im\Rm\ \lrar\ a^{\am +n}\equiv a^n\ (\mmod m). \]

\begin{thm} \label{rn20}
A number $a\in\Zm$ is regular iff the following logical inference holds
\[ p_i\mid a\ \rar\ \piai\mid a\ \ (i\in\N,\ \alpha_i>0). \]
\end{thm}
\noindent
\prf
If $a\in\Rm$ then $a\cdot a^{\phm}\equiv a\ (\mmod m)$. Let $i\in\N$ be such that $\alpha_i>0$, and suppose that $p_i\mid a$. Then as in the proof of Theorem \ref{in02}, we have
\[ a\equiv a\cdot a^{\phm}\equiv a\cdot 0\equiv 0\ (\mmod \piai). \]
Now, let us suppose that the inference holds. Let $i\in\N$ be such that $\alpha_i>0$. There are two possible cases.\\
1. $p_i\mid a$\\
Then we know that $\piai\mid a$ is true as well, so
\[ a\cdot a^{\phm}\equiv 0\equiv a\ (\mmod \piai). \]
2. $p_i\nmid a$\\
In this case we have $(a,\piai)=1$, so
\[ a\cdot a^{\phm}\equiv a\ (\mmod \piai). \]
From these two cases we have that $a\in\Rm$. \sqr

\begin{thm} \label{rn03}
A number $a\in\Zm$ is regular iff the following equivalence holds
\[ a^k\equiv a^l\ (\mmod m)\ \lrar\ k\equiv l\ (\mmod \am)\ \ (k,l\in\N). \]
It is also true that if $a\in\Zm$ is regular, then the following equivalence holds
\[ a^k\im\Em\ \lrar\ \am\mid k\ \ (k\in\N). \]
\end{thm}
\noindent
\prf
Let us first suppose that $a\in\Zm$ is regular. The $\rar$ part of the equivalence follows from Theorems \ref{nn02} and \ref{rn02}. The $\lar$ part is also true, since if $l\ge k$ and $k\equiv l\ (\mmod \am)$, then for some $q\ge 0$, we have $l=k+q\am$, so
\[ a^l\equiv a^{k+q\am}\equiv a^k a^{\am}\equiv a^k\ (\mmod m) \]
where the last congruence holds, because $a$ is regular.\\
Let us now suppose that the equivalence holds. Then, with $k:=\am +1,\ l:=1$, we have that $a$ is regular.\\
The second equivalence follows easily from the first. \sqr

\begin{thm} \label{rn21}
A number $a\in\Zm$ is regular iff there exists some $n>1$ such that
\[ a^n\equiv a\ (\mmod m). \]
\end{thm}
\noindent
\prf
First, let us suppose that such an $n$ exists. Then by Theorem \ref{in11} we have $a^{n-1}\im\Em$, from which $a^{n-1}\equiv a^{\am}\ (\mmod m)$ follows by the application of Theorem \ref{in03}. So, by multiplying this congruence by $a$, we get our desired result. Now, supposing that $a$ is regular, we may take $n:=\am +1$. \sqr

\begin{thm} \label{rn16}
For $a\in\Zm$
\[ a\in\Rm\ \lrar\ \left( a,\frac{m}{(a,m)} \right) =1\ \lrar\ (a,m)\in\Rm. \]
\end{thm}
\noindent
\prf
The first equivalence follows clearly from Theorem \ref{rn20}. Using this, and the fact that
\[ \left((a,m),\frac{m}{((a,m),m)}\right)=\left( a,\frac{m}{(a,m)} \right) \]
we get the second equivalence. \sqr

\begin{thm} \label{rn17}
For all $a\in\Nm$
\[ a\in\Rm\ \lrar\ (a^{-1})^{-1}\equiv a\ (\mmod m). \]
\end{thm}
\noindent
\prf
Follows from Theorem \ref{nn08}. \sqr

\begin{thm} \label{rn22}
Let $m_1,m_2\in\N$ be such that $m=[m_1,m_2]$, and $a\in\Zm$. Then
\[ a\in\Rm\ \lrar\ a\imo\Rmo\ \mand\ a\imt\Rmt. \]
Furthermore, if $a\in\Rm$ then $\am=[|a|_{m_1},|a|_{m_2}]$.
\end{thm}
\noindent
\prf
We will prove the equivalence using Theorem \ref{rn20}. Let us suppose that $m_1$ and $m_2$ have prime-factorizations
\[ m_1=\prod_{i=1}^{\infty} p_i^{\beta_i},\ \ m_2=\prod_{i=1}^{\infty} p_i^{\gamma_i}. \]
Then we have that $\alpha_i=\max(\beta_i,\gamma_i)$ for all $i\in\N$, since $m=[m_1,m_2]$.\\
The $\rar$ part: Taking any $i$ such that $p_i\mid a$, we have $\piai\mid a$, which implies $p_i^{\beta_i}\mid a$ and
$p_i^{\gamma_i}\mid a$, since $\alpha_i=\max(\beta_i,\gamma_i)$, so $a\imo\Rmo$ and $a\imt\Rmt$.\\
The $\lar$ part:
Taking any $i$ such that $p_i\mid a$, we have $p_i^{\beta_i}\mid a$ and $p_i^{\gamma_i}\mid a$,
which implies $\piai\mid a$, so we have that $a\in\Rm$.\\
Now, let us suppose that $a\in\Rm$. By Theorem \ref{in05}, for $k\in\N$ we have
\[ a^k\im\Em\ \lrar a^k\imo\Emo,\ a^k\imt\Emt\ \lrar\ [|a|_{m_1},|a|_{m_2}]\mid k. \]
From this we have that $[|a|_{m_1},|a|_{m_2}]\mid \am$. With $k=[|a|_{m_1},|a|_{m_2}]$ we have $a^{[|a|_{m_1},|a|_{m_2}]}\im\Em$ so $\am\le [|a|_{m_1},|a|_{m_2}]$, so $\am=[|a|_{m_1},|a|_{m_2}]$. \sqr

\begin{thm} \label{rn36}
For $e\in\Em$, there exist unique $m_1,m_2\in\N$ such that $m=m_1m_2$ and
\[ e\equiv 1\ (\mmod m_1)\ \mand\ e\equiv 0\ (\mmod m_2). \]
Furthermore, for all $a\in\Zm$, the following equivalence holds
\[ a\in\Rme\ \lrar\ (a,m_1)=1\ \mand\ m_2\mid a \]
and if $a\in\Rme$ then $\am = |a|_{m_1}$.
\end{thm}
\noindent
\prf
The proof is quite trivial, with the previous theorem in mind. \sqr

\begin{dfn} \label{rn37}
For $a\in\Zm,\ a^{\phm}\equiv e\ (\mmod m)$, using the notations of the theorem above, denote $\mma := m_1$.
\end{dfn}

We advise the reader to observe throughout our paper, that in many cases it is sufficient to examine the set $\mathrm{R}_{\mu_m(e)}^1$ instead of $\Rme$ (for any $e\in\Em$).

\begin{thm} \label{rn23}
For $a\in\Zm$
\[ a\in\Rm\ \lrar\ a\in_{\piai}\mathrm{R}_{\piai}\ \ (i\in\N) \]
and, if $a\in\Rm$ then
\[ \am=\mathrm{lcm} (|a|_{\piai}:1\le i\le\infty ). \]
Furthermore,
\[ |\Rm| = \prod_{\ali >0}(1+\varphi(\piai)). \]
\end{thm}
\noindent
\prf
Follows from Theorem \ref{rn22}. The formula for $|\Rm|$ follows from the fact that for $a\in\Z_{\piai}$
\[ a\in\R_{\piai}\ \lrar\ p_i\nmid a\ \mor\ \piai\mid a \]
which is a consequence of Theorem \ref{rn20}. \sqr

\begin{thm} \label{rn06}
For all $e\in\Em$, the structure $\langle \Rme; \{e,^{-1},\cdot\} \rangle$ is an Abelian group.
\end{thm}
\noindent
\prf
The properties to be shown are mostly trivial, except for maybe one. We need to show that for all $a\in\Rme$ there exists a unique $b\in\Rme$ such that $ab\equiv e\ (\mmod m)$. Let $b:=a^{-1}$. It is obvious that $b\in\Rme$ and $ab\equiv e\ (\mmod m)$. Now, let us suppose that there exists some other $b'\in\Rme$ such that $ab'\equiv e\ (\mmod m)$. Then we have
\[ a(b-b')\equiv 0\ (\mmod m)\ \rar\ 0\equiv a^{\am -1}\cdot a(b-b')\equiv \]
\[ \equiv e(b-b')\equiv b^{\bm +1}-(b')^{|b'|_m +1}\equiv b-b'\ (\mmod m).\ \ \sqr \]

\begin{thm} \label{rn42}
If $m$ is an odd number, then for all $e\in\Em$
\[ \prod \Rme \equiv (-1)^{2^{\omega(\mu_m(e))-1}}\cdot e\ (\mmod m). \]
\end{thm}
\noindent
\prf
We are going to make use of Theorem \ref{sd15} from the last section, which states that
\[ \prod \{a\in\Rme: \am\le 2\}\equiv (-1)^{2^{\omega(\mu_m(e))-1}}\cdot e\ (\mmod m). \]
For now, denote
\[ S:=\{a\in\Rme: \am\le 2\}. \]
Then for all $a\in\Rme$
\[ a\neq a^{-1}\ \lrar\ a\in\Rme\setminus S. \]
By our previous theorem, we have that
\[ \prod \Rme\setminus S \equiv e\ (\mmod m). \]
So we have
\[ \prod \Rme \equiv \left(\prod \Rme\setminus S\right)\cdot \left(\prod S\right) \equiv
e\cdot (-1)^{2^{\omega(\mu_m(e))-1}}\cdot e\ (\mmod m).\ \ \sqr \]

\vspace{0.25cm}
\begin{thm} \label{rn07}
For $a\in\Rm,\ n\in\N,\ i,j\in\Z$
\[ (a^n)^{-1} \equiv a^{-n}\ (\mmod m) \]
\[ a^{i+j} \equiv a^i\cdot a^j\ (\mmod m). \]
\end{thm}
\noindent
\prf
The first statement is equivalent to saying that
\[ a^{n\frac{\am}{(n,\am)}-n}\equiv a^{n\am -n}\ (\mmod m) \]
which by Theorem \ref{rn03} is equivalent to
\[ n\frac{\am}{(n,\am)}-n\equiv n\am -n\ (\mmod \am) \]
(when $n\frac{\am}{(n,\am)}-n\ne 0$), and this congruence obviously holds.\\
In the omitted case
\[ n\frac{\am}{(n,\am)}-n = 0\ \lrar\ \am\mid n. \]
So for some $k\in\N$, we have
\[ a^{n\am -n} = a^{(n-k)\am}\equiv a^0\ (\mmod \am). \]
For the second property, we can distinguish four different cases (for nonzero exponents).\\
Case of $i,j<0$:\\
\[ a^{i+j}=a^{-|i+j|}\equiv (a^{-1})^{|i+j|}=(a^{-1})^{|i|}\cdot (a^{-1})^{|j|}\equiv \]
\[ \equiv a^{-|i|}\cdot a^{-|j|}\equiv a^i\cdot a^j\ (\mmod m). \]
Case of $i<0,\ j>0$:\\
\textit{Case of $j\ge|i|$:}
\[ a^{i+j}=a^{j-|i|}\ \rar\ a^j=a^{i+j}\cdot a^{|i|}\ \rar \]
\[ \rar a^{i+j}\equiv a^j\cdot (a^{|i|})^{-1}\equiv a^j\cdot a^{-|i|}=a^i\cdot a^j\
(\mmod m). \]
\textit{Case of $j<|i|$:}
\[ a^{i+j}\equiv a^{j-|i|}\equiv a^{-(|i|-j)}\equiv (a^{|i|-j})^{-1}\equiv
(a^{|i|}\cdot a^{-j})^{-1}\ (\mmod m) \]
where the last congruence is true with the application of the previous case.\\
\[ (a^{|i|}\cdot a^{-j})\cdot (a^{-|i|}\cdot a^j)\equiv (a^{|i|}) (a^{|i|})^{-1}
(a^j)^{-1} (a^j)\equiv (a^{\am})^{|i|+j}\equiv a^{\am}\ (\mmod m). \]
So by the unicity of the inverse (Theorem \ref{rn06}), we have
\[ (a^{|i|}\cdot a^{-j})^{-1}\equiv a^{-|i|}\cdot a^j\equiv a^i\cdot a^j\ (\mmod m). \]
Case of $i,j>0$ is trivial.\\
Case of $i>0,\ j<0$ is similar to the case of $i<0,\ j>0$. \sqr

\begin{dfn} \label{rn08}
For $a\in\Zm$, let $\orbm{a}$ denote the set
\[ \{ a^n\ \mmod m:\ 1\le n\le \am\} \]
and in case of $A\subset\Zm$, let $\orbm{A}$ denote the set
\[ \bigcup_{a\in A} \orbm{a}. \]
\end{dfn}

\begin{thm} \label{rn09}
Let $b,c\in\Rm,\ n,k\in\N$. Then
\[ b^n,b^k\im\orbm{c}\ \lrar\ b^{(n,k)}\im\orbm{c}. \]
\end{thm}
\noindent
\prf
Let us first suppose that $b^n\equiv c^i,\ b^k\equiv c^j\ (\mmod m)$. Without hurting generality, we may suppose that there exist $x,y\ge 0$ such that $(n,k)=nx-ky$. So we have
\[ b^{(n,k)}=b^{nx-ky}=b^{nx+(-ky)}\equiv b^{nx}\cdot b^{-ky}\equiv b^{nx}\cdot
(b^{ky})^{-1}\equiv (c^{ix})\cdot (c^{jy})^{\phm -1}\im\orbm{c} \]
with the application of Theorem \ref{rn07}.\\
Now, let us suppose that $b^{(n,k)}\equiv c^l\ (\mmod m)$. Then we have
\[ b^n\equiv b^{(n,k)\frac{n}{(n,k)}}\equiv (c^l)^{\frac{n}{(n,k)}}\im\orbm{c} \]
The proof is similar for $b^k\im\orbm{c}$. \sqr

\begin{dfn} \label{rn10}
For $e\in\Em,\ b,c\in\Rme$, denote
\[ \D_m(b,c):= \mathrm{gcd} (n\in\N:\ 1\le n\le\bm,\ b^n\im\orbm{c}). \]
\end{dfn}

\begin{thm} \label{rn11}
If $e\in\Em,\ b,c\in\Rme$, then $\D_m(b,c)\mid\bm$ and
\[ b^k\im\orbm{c}\ \lrar\ \D_m(b,c)\mid k. \]
It is also true that $b^{\D_m(b,c)}\im\orbm{c}$ and
\[ \orbm{b}\cap\orbm{c} = \orbm{b^{\D_m(b,c)}}. \]
Furthermore
\[ |\orbm{b}\cap\orbm{c}|=\frac{\bm}{\D_m(b,c)}. \]
\end{thm}
\noindent
\prf
By Theorem \ref{rn09} and with induction, we have $b^{\D_m(b,c)}\im\orbm{c}$.\\
First, let us suppose that $\D_m(b,c)\mid k$. Then we have
\[ b^k\equiv (b^{\D_m(b,c)})^{\frac{k}{\D_m(b,c)}}\im\orbm{c}. \]
Now, if $b^k\im\orbm{c}$ then with $k':=k\ \mmod \bm$, we have $b^{k'}\im\orbm{c}$, so $\D_m(b,c)\mid k'$ by definition, and from this it follows that $\D_m(b,c)\mid k$.\\
So by the property now proven, we also have that
\[ \orbm{b}\cap\orbm{c} = \orbm{b^{\D_m(b,c)}}. \]
It is also true that $\D_m(b,c)\mid\bm$ since
\[ b^{\bm}\equiv e\equiv c^{\cm}\im\orbm{c}. \]
So we have
\[ |\orbm{b}\cap\orbm{c}|=|\orbm{b^{\D_m(b,c)}}|=|b^{\D_m(b,c)}|_m= \frac{\bm}{(\D_m(b,c),\bm)}=\frac{\bm}{\D_m(b,c)}.\ \ \sqr \]

\begin{thm} \label{rn13}
Let $k\in\N,\ e\in\Em,\ a,b\in\Rme,\ \eibma$. Then
\[ (k,\bm)\mid\ibma\ \lrar\ a^{\frac{\bm}{(k,\bm)}}\im\Em. \]
\end{thm}
\noindent
\prf
\[ e\equiv a^{\frac{\bm}{(k,\bm)}}\equiv b^{\frac{\bm\ibma}{(k,\bm)}}\ (\mmod m)\ \lrar \]
\[ \lrar\ \bm\mid \frac{\bm\ibma}{(k,\bm)}\ \lrar\ \frac{\ibma}{(k,\bm)}\in\Z\ \lrar\ (k,\bm)\mid\ibma.\ \ \sqr \]

\begin{thm} \label{rn14}
For $e\in\Em,\ a,b\in\Rme$
\[ (\am,\bm)=1\ \lrar\ |ab|_m=\am\cdot\bm \]
\[ \frac{[\am,\bm]}{(\am,\bm)}\mid |ab|_m\mid [\am,\bm] \]
\[ a^{\bm}\equiv b^{\am}\ (\mmod m)\ \rar\ \am=\bm. \]
\end{thm}
\noindent
\prf
We first prove the second property.
\[ (ab)^{[\am,\bm]}\equiv e\ (\mmod m)\ \rar\ |ab|_m\mid [\am,\bm] \]
\[ e\equiv (ab)^{\am\cdot |ab|_m}\equiv e\cdot b^{\am\cdot |ab|_m}
\equiv b^{\am\cdot |ab|_m}\ (\mmod m)\ \rar \]
\[ \rar\ \bm\mid \am\cdot |ab|_m\ \rar\ \frac{\bm}{(\am,\bm)}\mid |ab|_m \]
\[ e\equiv (ab)^{\bm\cdot |ab|_m}\equiv e\cdot a^{\bm\cdot |ab|_m}
\equiv a^{\bm\cdot |ab|_m}\ (\mmod m)\ \rar \]
\[ \rar\ \am\mid \bm\cdot |ab|_m\ \rar\ \frac{\am}{(\am,\bm)}\mid |ab|_m \]
\[ \rar\ \left[\frac{\am}{(\am,\bm)},\frac{\bm}{(\am,\bm)}\right]=\frac{[\am,\bm]}{(\am,\bm)}\mid |ab|_m. \]
The first property follows from the second one.\\
Now, we prove the third one.
\[ e\equiv a^{\bm\frac{\am}{(\am,\bm)}}\equiv b^{\frac{\am^2}{(\am,\bm)}}\ (\mmod m)\ \rar \]
\[ \rar\ \bm\mid \am\frac{\am}{(\am,\bm)}\ \lrar\ \frac{\bm}{(\am,\bm)}\mid\frac{\am}{(\am,\bm)}\ \rar\ \bm\mid\am \]
We get $\am\mid\bm$ the same way. \sqr

\begin{thm} \label{rn15}
Suppose that $e\in\Em,\ a,b,c\in\Rme$ and $a\in\orbm{b}\cap\orbm{c}$. Then there exists some $d\in\Rme$ for which $a\in\orbm{d}$ and $|d|_m=[\bm,\cm]$.
\end{thm}
\noindent
\textbf{Proof}\footnote{This proof was corrected in \cite{ba00135}.}\ \
By Theorem \ref{rn11}, we have
\[ \orbm{b}\cap\orbm{c}=\orbm{b^{\D_m(b,c)}}=\orbm{c^{\D_m(c,b)}} \]
so there exists some $K\in\N$ such that
\[ (b^{\D_m(b,c)})^K\equiv c^{\D_m(c,b)}\ (\mmod m) \]
and from
\[ |b^{\Db}|_m=|\orbm{b}\cap\orbm{c}|=|c^{\Dc}|_m=\frac{|b^{\Db}|_m}{(K,|b^{\Db}|_m)} \]
we have $(K,|b^{\Db}|_m)=1$.
\[ \frac{\bm}{\Db}=\frac{\cm}{\Dc}\ \rar\ \Dc\frac{\bm}{(\bm,\cm)}=\Db\frac{\cm}{(\bm,\cm)} \]
\[ \rar\ \frac{\bm}{(\bm,\cm)}\mid\Db\frac{\cm}{(\bm,\cm)} \]
and since
\[ \left(\frac{\bm}{(\bm,\cm )},\frac{\cm}{(\bm,\cm )}\right)=1 \]
we have
\[ \frac{\bm}{(\bm,\cm)}\mid\Db. \]
Denote
\[ l:=\frac{\Db (\bm,\cm )}{\bm} \]
\[ \bm=n_1 n_2 n_3,\ \cm=k_1 k_2 k_3,\ l=l_1 l_2 \]
where the numbers above have the properties
\[ k_1\mid n_1,\ n_2\mid k_2,\ l_1\mid n_2,\ l_2\mid k_1 \]
\[ 1=(n_3,k_3)=(n_i,n_j)=(n_i,k_j)=(k_i,k_j)\ \ (i\neq j). \]
Then we have
\[ (\bm,\cm)=n_2 k_1,\ |b^{n_2}|_m=n_1 n_3,\ |c^{k_1}|_m=k_2 k_3,\ (|b^{n_2}|_m,|c^{k_1}|_m)=1 \]
and also
\[ \Db=l\frac{n_1}{k_1} n_3,\ \Dc=l\frac{k_2}{n_2} k_3,\ |b^{\Db}|_m=\frac{n_2 k_1}{l}. \]
So by Theorem \ref{rn14}
\[ |b^{n_2} c^{k_1}|_m=n_1 n_3 k_2 k_3=[\bm,\cm ]. \]
Denote $d:=b^{n_2} c^{k_1}\ \mmod m$.
\[ d^{l\frac{n_1 k_2}{k_1 n_2} n_3 k_3}\equiv (b^{\Db})^{k_2 k_3} (c^{\Dc})^{n_1 n_3}\equiv (b^{\Db})^{k_2 k_3+Kn_1
n_3}\ (\mmod m) \]
\[ (k_2 k_3+Kn_1 n_3,\frac{n_2 k_1}{l_1 l_2})=1 \]
since $\frac{n_2}{l_1}\mid k_2$ but $(\frac{n_2}{l_1},Kn_1 n_3)=1$, and $\frac{k_1}{l_2}\mid n_1$ but $(\frac{k_1}{l_2},k_2 k_3)=1$. So there exists some $N\in\N$, such that
\[ (k_2 k_3+Kn_1 n_3)N\equiv 1\ (\mmod \frac{n_2 k_1}{l}). \]
So since $(b^{\Db})^I\equiv a\ (\mmod m)$ for some $I\in\N$, we have
\[ d^{l\frac{n_1 k_2}{k_1 n_2} n_3 k_3 NI}\equiv (b^{\Db})^{(k_2 k_3+Kn_1 n_3)NI}\equiv (b^{\Db})^I\equiv a\ (\mmod m).
\]
So $a\in\orbm{d}$. \sqr

\begin{thm} \label{rn41}
Let $e\in\Em,\ a\in\Rme,\ b\in\Nm$ be such that $a\in\orbm{b}$. Then there exists some $c\in\Rme$ such that $a\in\orbm{c}$ and
$\cm=\bm$.
\end{thm}
\noindent
\prf
Define $c:=be\ \mmod m$. Then obviously $c\in\Rme$, and
\[ c^{\ibma}\equiv b^{\ibma}\cdot e\equiv a\cdot e\equiv a\ (\mmod m). \]
So $a\in\orbm{c}$. Since $b^{\bm}\equiv e\ (\mmod m)$, we have $c\equiv b^{\bm +1}\ (\mmod m)$, so since $b\in\Nm$, we have
\[ \cm=\frac{\bm}{(\bm +1,\bm)}=\bm.\ \ \sqr \]

\begin{thm} \label{rn24}
Let $a\in\Rm,\ b\in\Nm$ be such that $\am\mid\bm$. Then the following inference holds
\[ \exists\ind_a^m\ b^n\ \rar\ \exists\ind_b^m\ a^n\ \ (n\in\N). \]
\end{thm}
\noindent
\prf
Let $n\in\N$ be such that $b^n\im\orbm{a}$. Then $\exists k\in\N:\ b^n\equiv a^k\ (\mmod m)$. From this we have
\[ |b^n|_m = |a^k|_m\ \rar\ (k,\am)\frac{\bm}{\am}=(n,\bm)\ \rar\ (k,\am)\mid n. \]
So $\exists l\in\N:\ kl\equiv n\ (\mmod \am)$, which implies $a^n\equiv a^{kl}\equiv b^{nl}\ (\mmod m)$. \sqr

It is easy to find numbers $m\in\N,\ a,b\in\Nm$ for which the theorem above does not hold.

\begin{thm} \label{rn25}
Let $a,b\in\Rm$ be such that $\am =\bm$. Then the following equivalence holds
\[ b^n\im\orbm{a}\ \lrar\ a^n\im\orbm{b}\ \ (n\in\N). \]
\end{thm}
\noindent
\prf
Follows easily from Theorem \ref{rn24}. \sqr

\begin{dfn} \label{rn39}
For $a,b\in\Rm$ we will say that $a$ and $b$ are equivalent modulo $m$ if the following are true
\[ a^{\am}\equiv b^{\bm}\ (\mmod m)\ \mand\ \am=\bm\ \mand\ \eibma \]
and we shall denote it as $a\eqm b$.
\end{dfn}

\begin{thm} \label{rn40}
The relation $\eqm$ truly is an equivalence relation.
\end{thm}
\noindent
\prf
It is clear that $a\eqm a$ for all $a\in\Rm$. By the previous theorem, we have that for all $a,b\in\Rm$
\[ a\eqm b\ \lrar\ \am=\bm,\ \eibma\ \lrar\ \am=\bm,\ \exists\ind_a^m b\ \lrar\ b\eqm a. \]
Now, to show the transitivity of our relation, we take any $a,b,c\in\Rm$
\[ a\eqm b,\ b\eqm c\ \rar\ \am=\bm=\cm,\ \exists n,k:\ b^n\equiv a,\ c^k\equiv b\ (\mmod m)\ \rar \]
\[ \rar\ \am=\cm,\ c^{nk}\equiv a\ (\mmod m)\ \rar\ a\eqm c.\ \ \sqr \]

\begin{thm} \label{rn26}
Let $a,b\in\Rm,\ n\in\N$ be such that $b^n\im\orbm{a},\ (n,\bm)=1$. Then $\orbm{b}\subset\orbm{a}$.
\end{thm}
\noindent
\prf
Let $k\in\N$ be such that $b^n\equiv a^k\ (\mmod m)$. Let $l$ be any integer. Then $\exists s\in\N:\ ns\equiv l\ (\mmod \bm)$. From this we have $b^l\equiv a^{ks}\ (\mmod m)$, which implies $b^l\im\orbm{a}$. \sqr

\begin{thm} \label{rn27}
Let $a,b\in\Rm,\ n\in\N$ be such that $a^n\im\orbm{b}$ and $\bm\mid\am,\ (n,\bm)=1$. Then $\orbm{b}\subset\orbm{a}$.
\end{thm}
\noindent
\prf
By Theorem \ref{rn24} we have that $b^n\im\orbm{a}$, so we have $\orbm{b}\subset\orbm{a}$ by Theorem \ref{rn26}. \sqr

\begin{thm} \label{rn28}
For $a,b\in\Rm$, we have $\orbm{a}=\orbm{b}$ iff $\am=\bm$ and there exists some $n\in\N$ such that $(n,\am)=1$ and
$a^n\im\orbm{b}$.
\end{thm}
\noindent
\prf
First, let us suppose that $\orbm{a}=\orbm{b}$. Then
\[ \am = |\orbm{a}| = |\orbm{b}| = \bm. \]
Also, for any $n\in\N,\ (n,\am)=1$, we have $a^n\im\orbm{a}=\orbm{b}$.\\
Now, let us suppose that $\am=\bm$ and there exists some $n\in\N$ such that $(n,\am)=1$ and $a^n\im\orbm{b}$. Then by Theorem \ref{rn26}, we have that $\orbm{a}\subset\orbm{b}$. By Theorem \ref{rn25}, we have that $b^n\im\orbm{a}$, so by applying Theorem \ref{rn26} once again, we get $\orbm{b}\subset\orbm{a}$ as well. \sqr

\begin{thm} \label{rn29}
Let $a,b\in\Rm$ be such that $\am\mid\bm$. Then there exists some $c\in\Rm$ such that $a,b\in\orbm{c}$, iff $a\in\orbm{b}$.
\end{thm}
\noindent
\prf
Let us suppose that $a,b\in\orbm{c}$ for some $c\in\Rm$. Then
\[ (\indm{c}{b},\cm)=\frac{\cm}{\bm}\mid \frac{\cm}{\am}=(\indm{c}{a},\cm)\mid \indm{c}{a}. \]
So there exists some $k\in\N$ such that $(\indm{c}{b})k\equiv \indm{c}{a}\ (\mmod \cm)$. So $b^k\equiv a\ (\mmod m)$. Now, if $a\in\orbm{b}$, then with $c:=b$, we have that $a,b\in\orbm{c}$. \sqr

\begin{thm} \label{rn32}
Let $e\in\Em,\ a,b\in\Rme$ be such that $(\am,\bm)=1$. Then $\orbm{a}\cap\orbm{b}=\{e\}$.
\end{thm}
\noindent
\prf
Let us take $c\in\orbm{a}\cap\orbm{b}$. Then there exists some $n\in\N$ such that $n\bm\equiv\indm{a}{c}\ (\mmod \am)$. So we have
\[ a^{n\bm^2}\equiv a^{\bm\indm{a}{c}}\equiv c^{\bm}\equiv (b^{\bm})^{\indm{b}{c}}\equiv e\ (\mmod m)\ \rar \am\mid n\bm^2\ \rar \]
\[ \rar \am\mid n\ \rar\ \indm{a}{c}\equiv n\bm\equiv 0\ (\mmod \am)\ \rar\ c\equiv a^{\indm{a}{c}}\equiv e\ (\mmod m).\ \ \sqr \]

\begin{thm} \label{rn30}
For $b\in\Rm,\ a\in\orbm{b},\ d\mid\bm$ the following equivalence holds
\[ \am =d\ \lrar\ \ibma = \frac{r\bm}{d},\ (r,d)=1. \]
\end{thm}
\noindent
\prf
First, let us suppose that $\am =d$. Then
\[ d=\am =\frac{\bm}{(\ibma,\bm)}\ \rar\ \ibma =\frac{r\bm}{d},\ r:=\frac{\ibma}{(\ibma,\bm)} \]
where we see that $(r,\am)=1$.\\
Now, let $r,d,b\in\N$ be any numbers such that $(r,d)=1,\ d\mid b$. Then
\[ \frac{b}{\left(\frac{rb}{d},b\right)} =d. \]
Considering the relation between prime-factorizations and the greatest common divisor, we need to show that for all $r,d,b\ge 0,\ \min(r,d)=0,\ d\le b$, we have $b-\min(r+b-d,b)=d$. Supposing that $r\le d$, we get $r=0$, from which the desired relation is $b-(b-d)=d$. Now, if $d\le r$, then $d=0$, so the relation we need is $b-b=0$.\\
Applying the relation above, we get that
\[ \am =\frac{\bm}{\left(\frac{r\bm}{d},\bm\right)} =d.\ \ \sqr \]

\begin{thm} \label{rn38}
For $a\in\Rm,\ d\mid\am$
\[ |\{b\in\orbm{a}:\ \bm=d\}| = \varphi(d) \]
and
\[ |\{b\in\Rm:\ b\eqm a\}| = \varphi(\am). \]
\end{thm}
\noindent
\prf
The first equality follows from the above theorem. The second one follows from the first with $d=\am$. \sqr

We now examine further properties of the function $\D_m(a,b)$.

\begin{dfn} \label{rn34}
For $e\in\Em,\ a,b\in\Rme$, let the relative order of the numbers $a$ and $b$ be defined as
\[ |a,b|_m := |\orbm{a}\cap\orbm{b}| = \frac{\am}{\D_m(a,b)} = \frac{\bm}{\D_m(b,a)}. \]
\end{dfn}

\begin{thm} \label{rn33}
For $e\in\Em,\ a,b\in\Rme,\ n\in\N$
\[ \D_m(a,b)=\D_m(b,a)\ \lrar\ \am=\bm \]
\[ \frac{\D_m(a,b)}{\D_m(b,a)} = \frac{\am}{\bm} \]
\[ \D_m(a,b) = \left| a^{|a^{\D_m(a,b)}|_m} \right|_m \]
\[ \D_m(a^n,b) = \frac{\D_m(a,b)}{(n,\D_m(a,b))} \]
\[ (|a,b|_m,\D_m(b,a))=1\ \rar\ \D_m(a,b^n) = \D_m(a,b)\left( n,\frac{\am}{\D_m(a,b)} \right). \]
\end{thm}
\noindent
\prf
The first, second, and third relations follow trivially from Theorem \ref{rn11}. To prove the fourth relation, for any $k\in\N$ we see that
\[ \D_m(a^n,b)\mid k\ \lrar\ a^{nk}\im\orbm{b}\ \lrar\ \D_m(a,b)\mid nk\ \lrar\ \frac{\D_m(a,b)}{(n,\D_m(a,b))}\mid k. \]
To prove the fifth relation,
\[ \D_m(a,b^n) = \frac{\am}{|b^n|_m} \D_m(b^n,a) = \frac{\am}{\bm} (n,\bm)\frac{\D_m(b,a)}{(n,\D_m(b,a))} = \]
\[ = \D_m(a,b)\frac{(n,\bm)}{(n,\D_m(b,a))} = \D_m(a,b)\left( n,\frac{\bm}{\D_m(b,a)}\right) = \]
\[ = \D_m(a,b)\left(n,\frac{\am}{\D_m(a,b)}\right).\ \ \sqr \]

\begin{thm} \label{rn35}
For $e\in\Em,\ a,b\in\Rme,\ n,k\in\N$
\[ |a,b|_m = |b,a|_m \]
\[ |a,a|_m = \am \]
\[ (\am,\bm)=1\ \rar\ |a,b|_m =1 \]
\[ b\in\orbm{a}\ \rar\ |a,b|_m = \bm \]
\[ (|a,b|_m,\D_m(a,b))=1\ \rar\ |a^n,b|_m = \frac{|a,b|_m}{(n,|a,b|_m)} \]
\[ (|a,b|_m,\D_m(a,b^k)\D_m(b,a))=1\ \rar\ |a^n,b^k|_m = \frac{|a,b|_m}{(n(k,|a,b|_m),|a,b|_m)}. \]
\end{thm}
\noindent
\prf
The first and second relations are trivial. The first inference follows from Theorem \ref{rn32}, and the second is trivial. Now, to prove the fifth relation,
\[ |a^n,b|_m = \frac{\bm}{\D_m(b,a^n)} = \frac{\bm}{\D_m(b,a)\left(n,\frac{\bm}{\D_m(b,a)}\right)} = \frac{|a,b|_m}{(n,|a,b|_m)}. \]
Lastly, since
\[ (|a,b|_m,D_m(a,b^k))=1\ \rar\ (|a,b^k|_m,D_m(a,b^k))=1 \]
we have
\[ |a^n,b^k|_m = \frac{|a,b^k|_m}{(n,|a,b^k|_m)} = \frac{|b,a|_m}{(k,|b,a|_m)\left(n,\frac{|b,a|_m}{(k,|b,a|_m)}\right)} = \]
\[ = \frac{|a,b|_m}{(n(k,|a,b|_m),|a,b|_m)}. \ \ \sqr \]

\newpage
\section{Binomial Congruences}

\begin{dfn} \label{bc02}
For $a\in\Z,\ k\in\N$ let $\Mmka$ denote the logical function which gives ``true'' if the equation in $x$
\[ x^k\equiv a\ (\mmod m) \]
is solvable, otherwise let its value be ``false''. Furthermore, let $\Sm(k,a)$ denote the set of solutions in $\Zm$ of the equation above, and let $\SmRka$ denote the set of regular solutions in \Zm. In case of $a\in\Rm$, let \oma\ denote the number
\[ \max\ \{\bm:\ b\in\Rm,\ \eibma\} \]
furthermore let $\ind^m a:=\frac{\oma}{\am}$.
\end{dfn}

Note that $\oma =\phm$ for some $a\in\Rm$, iff $\Rm^1$ is cyclical.

\begin{thm} \label{bc01}
For $a\in\Rm,\ k\in\N$
\[ \Mmka\ \lrar\ a^{\frac{\oma}{(k,\oma)}}\im\Em. \]
\end{thm}
\noindent
\prf
Let $b\in\Rm$ be such that $\eibma$ and $\bm=\oma$. Then we have
\[ a^{\frac{\oma}{(k,\oma)}}\im\Em\ \lrar\ (k,\bm)\mid\ibma. \]
If $(k,\bm)\mid\ibma$ holds, then there must exist some $1\le l\le\bm$ for which $kl\equiv \ibma\ (\mmod \bm)$. So we have
\[ b^{kl}\equiv b^{\ibma}\ (\mmod m)\ \rar\ (b^l)^k\equiv a\ (\mmod m) \]
so $b^l$ is a solution of the equation.\\
Now, let $x_0$ be a solution of the equation, and denote $e:= a^{\phm}\ \mmod m,\ \ c:=x_0 e\ \mmod m$. Then we have that $c$ is a solution as well, since
\[ c^k\equiv (x_0)^k e\equiv a\cdot a^{\am}\equiv a\ (\mmod m) \]
and $c\in\Rm$ since
\[ c\cdot c^{\cm}\equiv c\cdot c^{\phm}\equiv x_0 e (x_0^k)^{\phm}\equiv x_0 e\equiv c\ (\mmod m). \]
It is also clear that $\cm\mid\bm$. For, let us make the indirect assumption that $\cm\nmid\bm$. Then we have $\cm < \bm$ by the definition of \oma. We also know by Theorem \ref{rn15} that there exists some $d\in\Rm$, such that $\exists\ind_d^m a$ and $|d|_m=[\bm,\cm]$. It is clear that $|d|_m>\bm$, which obviously contradicts the selection of $b$ and the definition of \oma.
So we must have that $\cm\mid\bm$. From this, we have
\[ a^{\frac{\oma}{(k,\oma)}}\equiv a^{\frac{\bm}{(k,\bm)}}\equiv (c^k)^{\frac{\bm}{(k,\bm)}}\equiv (c^{\cm})^{\frac{\bm}{\cm}\cdot\frac{k}{(k,\bm)}}\equiv e\ (\mmod m).\ \ \sqr \]

The difficulty of the verification of the condition
\[ a^{\frac{\oma}{(k,\oma)}}\im\Em \]
lies within the calculation of $\oma$. So I believe that the examination of the mapping $m\mapsto \oma$ is probably the most logical direction, research on this subject should take. Let us look at some immediate corollaries of our theorem.

\begin{thm} \label{bc09}
For $a\in\Rm,\ k\in\N$
\[ \Mmka\ \lrar\ \M_m((k,\phm),a)\ \lrar\ \M_m((k,\psm),a). \]
\end{thm}
\noindent
\prf
The equivalence follows trivially from our previous theorem, since
\[ \oma\mid\phm,\psm\ \rar\ (k,\oma)=((k,\phm),\oma)=((k,\psm),\oma).\ \ \sqr \]

\begin{thm} \label{bc08}
For $a\in\Rm,\ k_1,k_2\in\N$
\[ \M_m(k_1,a)\ \mand\ \M_m(k_2,a)\ \lrar\ \M_m([k_1,k_2],a). \]
\end{thm}
\noindent
\prf
Our theorem follows from Theorem \ref{rn09}, and the fact that
\[ \left(\frac{\oma}{(k_1,\oma)},\frac{\oma}{(k_2,\oma)}\right) = \frac{\oma}{([k_1,k_2],\oma)}.\ \ \sqr \]

We now look at a necessary and then a sufficient condition for the solvability of a binomial congruence modulo $m$.

\begin{thm} \label{bc03}
For $a\in\Zm,\ k\in\N$
\[ \Mmka\ \rar\ a^{\frac{\phm}{(k,\phm)}}\im\Em. \]
\end{thm}
\noindent
\prf
Let the solution of the binomial congruence, be denoted by $x_0$. Then
\[ a^{\frac{\phm}{(k,\phm)}} \equiv (x_0^k)^{\frac{\phm}{(k,\phm)}} \equiv (x_0^{\phm})^{\frac{k}{(k,\phm)}} \im\Em.\ \ \sqr \]

\begin{thm} \label{bc04}
Let $a,b\in\Rm,\ k\in\N$ be such that $\eibma$ and $(k,\bm)\mid\ibma$. Then \Mmka.
\end{thm}
\noindent
\prf
If the conditions above are satisfied, then for some $l\in\Z_{\bm}$, we have $kl\equiv \ibma\ (\mmod \bm)$. So since $b\in\Rm$, we have
\[ b^{kl}\equiv b^{\ibma}\ (\mmod m)\ \rar\ b^l\im\Sm(k,a)\ \rar\ \Mmka.\ \ \sqr \]

We now look at some special solutions of a binomial congruence.

\begin{thm} \label{bc05}
Let $a,b\in\Rm,\ e\in\E_{\bm},\ k,l\in\N,\ kl\in\R_{\bm}^e$ be such that $\eibma\in\R_{\bm}^e,\ kl\equiv e\ (\mmod
\bm)$ and \Mmka. Then
\[ b^{l\ibma +n\frac{\bm}{(k,\bm)}}\im\Sm(k,a)\ \ (n\in\N). \]
\end{thm}
\noindent
\prf
\[ k(l\ibma +n\frac{\bm}{(k,\bm)})\equiv \ibma + n\frac{k}{(k,\bm)}\bm \equiv \ibma\ (\mmod \bm)\ \rar \]
\[ \rar (b^{l\ibma +n\frac{\bm}{(k,\bm)}})^k\equiv b^{\ibma}\equiv a\ (\mmod m).\ \ \sqr \]

Next, we examine the number of solutions of a binomial congruence.

\begin{thm} \label{bc07}
If $a\in\Rm,\ k\in\N$ and \Mmka, then $|\SmRka|>0$.
\end{thm}
\noindent
\prf
It is clear that $x_0\cdot e\im\SmRka$ for any $x_0\in\Sm(k,a)$, where $a^{\am}\equiv e\ (\mmod m)$. \sqr

\begin{thm} \label{bc06}
If $e\in\Em,\ a\in\Rme,\ k\in\N$ and \Mmka, then $|\SmRka|=|\SmR(k,e)|$.
\end{thm}
\noindent
\prf
Let $x_0\in\SmRka$ be some regular solution. Then according to Theorem \ref{rn06}, we have exactly one $\xni\in\Rme$ such that $\xni x_0\equiv e\ (\mmod m)$. Let us define the set
\[ A:=\{\xni x_i:\ x_i\in\SmRka \}. \]
Then we have that $A\subset\SmR(k,e)$, since for any $x_i\in\SmRka$
\[ (\xni x_i)^k \equiv (\xni)^k x_i^k \equiv (\xni)^k x_0^k \equiv e\ (\mmod m) \]
and for $i\neq j$ we have $\xni x_i\ \not\equiv\ \xni x_j\ (\mmod m)$, for let suppose that for some $i\neq j$
\[ \xni x_i\equiv \xni x_j\ \rar\ x_i\equiv x_i \xni x_0\equiv x_j \xni x_0 \equiv x_j\ (\mmod m) \]
which is a contradiction. So we have that $|\SmRka|=|A|\le |\SmR(k,e)|$.\\
We also have that
\[ \SmRka = x_0\cdot \SmR(k,e)\ \mmod m \subset \SmRka \]
so $|\SmR(k,e)|\le |\SmRka|$. \sqr

\newpage
\section{Generalized Primitive Roots}

\begin{dfn} \label{pr01}
A number $g\in\Rm$ is said to be a generalized primitive root modulo $m$, if $\omega_m(g)=|g|_m$. Let the set of such $g$ be denoted by $\Gm$. Furthermore, let $\Oma$ denote the set
\[ \{b\in\Rm:\ \eibma\ \mand\ \bm=\oma\}. \]
\end{dfn}

\begin{thm} \label{pr02}
For all $a\in\Rm$
\[ \Oma\subset\Gm. \]
Furthermore, for all $g\in\Gm$ there exists some $a\in\Rm$ such that $g\in\Oma$.
\end{thm}
\noindent
\prf
To prove the first part of our theorem, take any $b\in\Oma$, and let us suppose indirectly that $b\notin\Gm$. Then there exists some $c\in\Rm$, such that $\exists\indm{c}{b}$ and $\cm>\bm$. Obviously $\exists\indm{c}{a}$, since
\[ (c^{\indm{c}{b}})^{\indm{b}{a}}\equiv b^{\indm{b}{a}}\equiv a\ (\mmod m). \]
Which contradicts the maximality of $\bm=\oma$.\\
Now, to prove the second part, take any $g\in\Gm$. It is trivial, that $g\in\Omega_m(g)$. \sqr

Note that our theorem implies the nonemptyness of $\Gm$.

\begin{thm} \label{pr05}
For $a\in\Rm,\ g\in\Oma$ the following equivalence holds
\[ g^n\im\Oma\ \lrar\ (n,|g|_m)=1\ \ (n\in\N). \]
\end{thm}
\noindent
\prf
To prove the $\rar$ part of the equivalence
\[ |g|_m = \oma = |g^n|_m = \frac{|g|_m}{(n,|g|_m)}\ \rar\ (n,|g|_m)=1. \]
Now, if we suppose that $(n,|g|_m)=1$, then there exists some $k\in\N$ such that
$nk\equiv\indm{g}{a}\ (\mmod |g|_m)$. So
\[ (g^n)^k\equiv g^{\indm{g}{a}}\equiv a\ (\mmod m)\ \rar\ \exists\indm{g^n}{a} \]
and
\[ |g^n|_m=\frac{|g|_m}{(n,|g|_m)}=|g|_m=\oma. \]
So $g^n\im\Oma$. \sqr

\begin{thm} \label{pr06}
For $a\in\Rm$ the following equivalence holds
\[ g\in\Oma\ \lrar\ g^{-1}\in\Oma. \]
\end{thm}
\noindent
\prf
Follows from our previous theorem. \sqr

\begin{thm} \label{pr04}
Let $m_1,m_2\in\N$ be such that $m=[m_1,m_2]$, and $g\in\Zm$. If $g\imo\Gmo$ and $g\imt\Gmt$, then $g\im\Gm$.
\end{thm}
\noindent
\prf
Let us suppose indirectly, that $g\notin\Gm$. This means that there exists some $h\in\Rm$ and $n\in\N$, such that $|h|_m>|g|_m$, and $h^n\equiv g\ (\mmod m)$. This implies that $(n,|h|_m)>1$, and for $i\in\{1,2\}$
\[ h^n\equiv g\ (\mmod m_i)\ \rar\ |g|_{m_i}=\frac{|h|_{m_i}}{(n,|h|_{m_i})}\ \rar\ (n,|h|_{m_i})=1. \]
So combining the two we get
\[ 1=[(n,|h|_{m_1}),(n,|h|_{m_2})]=(n,[|h|_{m_1},|h|_{m_2}])=(n,|h|_m) \]
which is a contradiction. \sqr

Our theorem above sheds some light on the still hazy structure of $G_m$.

\begin{thm} \label{pr03}
If $a,b\in\Rm,\ a\eqm b$, then $\oma=\omega_m(b)$.
\end{thm}
\noindent
\prf
It is clear, that $(\ibma,\bm)=1$. So for any $g\in\Omega_m(b)$
\[ (g^{\indm{g}{b}})^{\ibma}\equiv a\ (\mmod m) \]
so $\exists\indm{g}{a}$. So
\[ \omega_m(b)=|g|_m\le \oma. \]
The inequality $\oma\le\omega_m(b)$ may be proven in the same way. \sqr

The theorem above shows, that in our quest of finding an easy method for the calculation of the function $\omega_m$, it would
be worth examining the equivalence classes according to the relation $\eqm$. It also implies that if a number is equivalent to a gen. primitive root, then it is a gen. primitive root as well. Therefore, it would also be worth examining the structure of $\Gm$ and $\Oma$, partitioned according to our equivalence relation.

\newpage
\section{Number Theoretic Functions}

In this section, we shall examine some functions, along with some of their properties, emerging from
the discussions above.

\begin{dfn} \label{fs08}
A number theoretic function $f:\N\to\N$ is said to be multiplicative, if for all $a,b\in\Domf,\ (a,b)=1$
\[ a\cdot b\in\Domf\ \mand\ f(a\cdot b) = f(a)\cdot f(b) \]
and we shall denote it as $f\in\MU$. We will say that $f$ is quasimultiplicative, if for all $a,b\in\Domf$
\[ [a,b]\in\Domf\ \mand\ f([a,b]) = [f(a),f(b)] \]
and we shall denote it as $f\in\QM$. We will say that $f$ is division-invariant, if the following inference holds for all $a,b\in\Domf$
\[ a\mid b\ \rar\ f(a)\mid f(b) \]
and we shall denote it as $f\in\DI$. Lastly, we will say that $f$ is prime-power division-invariant, if for all primes $q$ and $\beta,\gamma\in\N,\ \beta\le\gamma$, such that $q^{\beta},q^{\gamma}\in\Domf$, we have $f(q^{\beta})\mid f(q^{\gamma})$; and we shall denote it as $f\in\DIpa$.
\end{dfn}

Note that the functions $\psi$ (by our theorem below), $m\mapsto \am$ (with domain $\{m\in\N: a\im\Nm\}$), $a\mapsto (a,b)$ are quasimultiplicative. We suspect, that for most (if not all) quasimultiplicative functions, there exists some quick algorithm for their computation. The basis of this conjecture is that the well-known Euclidean Algorithm computes the function  $a\mapsto (a,b)\in\QM$. Furthermore, it is also possible, that the computation of most multiplicative functions relies heavily on prime-factorization; that is, their computation is mostly equivalent to prime-factorization, in terms of speed.

\begin{thm} \label{fs13}
For any $g\in\DIpa,\ \Dom(g)=\N$ and $n\in\N$, with prime-factorization $n=\prod_{i\in\N} p_i^{\gamma_i}$, define the function $f$ as
\[ f(n) = \lcm(g(p_i^{\gamma_i}):\ i\in\N). \]
Then $f\in\QM$.
\end{thm}
\noindent
\prf
Take $a,b\in\Domf$, with prime-factorizations $a=\prod_{i\in\N} p_i^{\gamma_i},\ b=\prod_{i\in\N} p_i^{\delta_i}$. Then
\[ f([a,b]) = \lcm(g(p_i^{\max(\gamma_i,\delta_i)}):\ i\in\N) = \]
\[ = \lcm([g(p_i^{\min(\gamma_i,\delta_i)}), g(p_i^{\max(\gamma_i,\delta_i)})]:\ i\in\N) = \]
\[ = [ \lcm(g(p_i^{\gamma_i}):\ i\in\N), \lcm(g(p_i^{\delta_i}):\ i\in\N) ] = [f(a),f(b)].\ \ \sqr \]

It is interesting to ponder the question whether there would exist such a $g$ for all $f\in\QM$.

\begin{thm} \label{fs09}
\[ f\in\QM\ \lrar\ f\in\DI\ \mand\ (a,b\in\Domf,\ (a,b)=1\ \rar\ f(ab)=[f(a),f(b)]). \]
\end{thm}
\noindent
\prf
First, let us suppose that $f\in\QM$, and take $a,b\in\Domf$ such that $a\mid b$. Then
\[ f(b) = f([a,b]) = [f(a),f(b)]\ \rar\ f(a)\mid f(b). \]
Now, suppose that the right hand side of the equivalence holds. Let $a,b\in\Domf$ and $n_i,k_i\in\N\ (i=1,2,3)$ be such that
\[ a=n_1 n_2 n_3,\ b=k_1 k_2 k_3,\ n_1\mid k_1,\ k_2\mid n_2 \]
\[ 1=(n_i,n_j)=(k_i,k_j)=(n_i,k_j)=(n_3,k_3)\ (i\neq j). \]
Such decompositions exist, and are easy to find, by looking at the prime factorizations of $a$ and $b$. So we have $[a,b]=k_1 n_2 n_3 k_3$, and
\[ f([a,b]) = [f(k_1),f(n_2),f(n_3),f(k_3)] = [[f(n_1),f(k_1)],[f(k_2),f(n_2)],f(n_3),f(k_3)] = \]
\[ = [[f(n_1),f(n_2),f(n_3)],[f(k_1),f(k_2),f(k_3)]] = [f(a),f(b)].\ \ \sqr \]

\begin{thm} \label{fs10}
\[ f\in\DI\ \lrar\ \forall a,b\in\Domf:\ [f(a),f(b)]\mid f([a,b]). \]
\end{thm}
\noindent
\prf
First, let us suppose that $f\in\DI$. For any $a,b\in\Domf$
\[ a,b\mid [a,b]\ \rar\ f(a),f(b)\mid f([a,b])\ \rar\ [f(a),f(b)]\mid f([a,b]). \]
Now, suppose that the right hand side property is what holds for $f$. Then for any $a,b\in\Domf$
\[ a\mid b\ \rar\ f([a,b])=f(b)\ \rar\ [f(a),f(b)]\mid f(b)\ \rar\ [f(a),f(b)]=f(b)\ \rar\ f(a)\mid f(b).\ \ \sqr \]

\begin{thm} \label{fs11}
If $f\in\QM$ is injective, then the following equivalence holds
\[ a\mid b\ \lrar\ f(a)\mid f(b)\ \ (a,b\in\Domf). \]
\end{thm}
\noindent
\prf
If $f\in\QM$, then by Theorem \ref{fs09}, we have the $\rar$ part of the equivalence. Now, suppose that $a,b\in\Domf$ and $f(a)\mid f(b)$. Then
\[ f(b) = [f(a),f(b)] = f([a,b])\ \rar\ b=[a,b]\ \rar\ a\mid b.\ \ \sqr \]

\begin{dfn} \label{fs01}
For $e\in\Em,\ k\in\N$, let us define the following sets
\[ \kRm:=\{a\in\Rm:\ \am =k \},\ \ \kRme:=\ \kRm\cap\Rme. \]
Now, define
\[ r_m^e(k):= |\kRme|\ \ (k\in\N). \]
Furthermore, let
\[ \rhme (k):= |\{a\in\Rme:\ \am\mid k \}|\ \ (k\in\N). \]
\end{dfn}

Note that by Theorem \ref{bc06} we have that
\[ |\SmRka| = |\SmR(k,e)| = \rhme(k) \]
if $k\in\N,\ e\in\Em,\ a\in\Rme$ and $\Mmka$.

\begin{thm} \label{fs02}
For all $e\in\Em$
\[ r_m^e(k)=r_{\mu_m(e)}^1(k),\ \ \rhme(k)=\rho_{\mu_m(e)}^1(k)\ \ (k\in\N). \]
Furthermore, if $m$ is weakly even, then $r_m^e,\rhme\in\MU$.
\end{thm}
\noindent
\prf
By the application of Theorem \ref{rn36}, we have the first statement of the theorem. So, this result shows that it is enough to prove the multiplicativity of our functions for the case of $e=1$.\\
Let $k_1,k_2\in\N$ be such that $(k_1,k_2)=1$ and $r_m^1(k_1),r_m^1(k_2)>0$ (otherwise the theorem holds trivially). Let $k:=k_1 k_2$ and $n:=\om$. The following equivalence is quite trivial. For all $x\in\N^n$
\[ [x]=k\ \lrar\ \exists! u,v\in\N^n:\ x_i=u_i v_i,\ (u_i,v_i)=1\ (1\le i\le n)\ \mand\ k_1=[u],\ k_2=[v]. \]
For simplicity's sake, let us suppose that $\alpha_i >0$ for all $1\le i\le n$. Let $f$ denote the integral vector $(\varphi(p_1^{\alpha_1}),\dots,\varphi(p_n^{\alpha_n}))$.\\
Since $m$ is weakly even, we know that for all moduli $\piai\ (1\le i\le n)$, there exists a primitive root modulo $\piai$. So, it is well-known, that if a primitive root exists modulo \piai, then the number of integers in $\Z_{\piai}$ of order
$d\in\N$ is $\varphi(d)$ (for $d\mid\varphi(\piai)$). Now, with the above facts, and the Chinese Remainder Theorem in mind, we have the following
\[ r_m^1(k) = \sum \left(\prod_{i=1}^n \varphi(x_i):\ x\in\N^n,\ [x]=k,\ x\mid f \right) = \]
\[ = \sum \left(\prod_{i=1}^n \varphi(x_i):\ u\in\N^n,\ [u]=k_1,\ u\mid f,\ v\in\N^n,\ [v]=k_2,\ v\mid f \right) = \]
\[ = \sum \left(\left(\prod_{i=1}^n \varphi(u_i)\right) \left(\prod_{i=1}^n \varphi(v_i)\right):
\ u\in\N^n,\ [u]=k_1,\ u\mid f,\ v\in\N^n,\ [v]=k_2,\ v\mid f \right) = \]
\[ = \left(\sum \left(\prod_{i=1}^n \varphi(u_i):\ u\in\N^n,\ [u]=k_1,\ u\mid f\right)\right) \cdot \]
\[ \cdot \left(\sum \left(\prod_{j=1}^n \varphi(v_j):\ v\in\N^n,\ [v]=k_2,\ v\mid f\right)\right) = \]
\[ = r_m^1(k_1)\cdot r_m^1(k_2). \]
The multiplicativity of the function $\rho_m^1$, may be shown similarly. All we need to do, is change some equality signs to division signs, as follows,
\[ \rhmo(k) = \sum \left(\prod_{i=1}^n \varphi(x_i):\ x\in\N^n,\ [x]\mid k,\ x\mid f \right) = \]
\[ = \sum \left(\prod_{i=1}^n \varphi(x_i):\ u\in\N^n,\ [u]\mid k_1,\ u\mid f,\ v\in\N^n,\ [v]\mid k_2,\ v\mid f \right) = \]
\[ = \sum \left(\left(\prod_{i=1}^n \varphi(u_i)\right) \left(\prod_{i=1}^n \varphi(v_i)\right):
\ u\in\N^n,\ [u]\mid k_1,\ u\mid f,\ v\in\N^n,\ [v]\mid k_2,\ v\mid f \right) = \]
\[ = \left(\sum \left(\prod_{i=1}^n \varphi(u_i):\ u\in\N^n,\ [u]\mid k_1,\ u\mid f\right)\right) \cdot \]
\[ \cdot \left(\sum \left(\prod_{j=1}^n \varphi(v_j):\ v\in\N^n,\ [v]\mid k_2,\ v\mid f\right)\right) = \]
\[ = r_m^1(k_1)\cdot r_m^1(k_2).\ \ \sqr \]

So the theorem above, tells us, that if $m$ is weakly even, then it is enough to determine $r_m^e$ and $\rhme$ at the prime-power divisors of $k\in\N$.

\begin{thm} \label{fs03}
Supposing that $m$ is weakly even, $\beta\in\N$, $q\in\N$ is prime, $q^{\beta}\mid\psm,\ n=\om$, and without hurting generality, we may also suppose that for some $\delta\in (\N\cup\{0\})^n$
\[ \alpha_i>0,\ \varphi(\piai)=q^{\delta_i} r_i,\ (q,r_i)=1\ \ (1\le i\le n). \]
Then
\[ \rhmo (q^{\beta}) = q^{\sum_{i=1}^n \min(\beta,\delta_i)} \]
\[ \rmo (q^{\beta}) = \rhmo (q^{\beta}) - \rhmo (q^{\beta -1}) \]
\[ \rmo (q) = q^{\Delta} -1,\ \ \Delta=|\{\delta_i\neq 0:\ 1\le i\le n\}|. \]
\end{thm}
\noindent
\prf
The first relation may be proven in the following manner, by considering the ideas that follow from $m$ being weakly even, like in the proof of the previous theorem,
\[ \rhmo (q^{\beta}) = \sum \left( \prod_{i=1}^n \varphi(q^{\gamma_i}):\ \max(\gamma)\le \beta,\ \gamma\le\delta \right) = \]
\[ = \sum_{\gamma_i\le \min(\beta,\delta_i)} q^{\sum_i \gamma_i} \cdot \left(1-\frac{1}{q}\right)^{\sum_{\gamma_i\neq 0} 1} = \]
\[ = \prod_{i=1}^n \left( 1 + q\left(1-\frac{1}{q}\right) + q^2\left(1-\frac{1}{q}\right) + \dots + q^{\min(\beta,\delta_i)}
\left(1-\frac{1}{q}\right) \right) = \]
\[ = \prod_{i=1}^n \left( 1 + \frac{q-1}{q}\left(-1+\frac{q^{\min(\beta,\delta_i)+1} -1}{q-1}\right) \right) = \]
\[ = \prod_{i=1}^n \left( 1-1+\frac{1}{q}+ q^{\min(\beta,\delta_i)} -\frac{1}{q} \right) = q^{\sum_{i=1}^n \min(\beta,\delta_i)}. \]
The second relation is quite trivial. The third one follows from the first and the second. \sqr

\begin{thm} \label{fs04}
If $m$ is weakly even, then
\[ \rhmo(k) = \prod_{i=1}^{\infty} (k,\varphi(\piai))\ \ (k\in\N). \]
\end{thm}
\noindent
\prf
First, let us examine the case of $m=p^{\alpha},\ \alpha\in\N,\ p$ is prime, and there exists a primitive root modulo $p^{\alpha}$, and $k=q^{\beta},\ q\ prime,\ \beta\in\N$.\\
Then, for some $r\in\N$, we have $\varphi(p^{\alpha})=q^{\delta} r,\ q\nmid r$, so by our previous theorem, we get
\[ \rho_{p^{\alpha}}^1(q^{\beta}) = q^{\min(\beta,\delta)} = (q^{\beta},\varphi(p^{\beta})). \]
The general case follows quite trivially, through the Chinese Remainder Theorem. \sqr

\begin{thm} \label{fs12}
If $m$ is weakly even, then $\rhme\in\DI$ for all $e\in\Em$.
\end{thm}
\noindent
\prf
Follows easily from our previous theorem. \sqr

\begin{thm} \label{fs05}
For $e\in\Em,\ k\in\N$
\[ |\orbm{\kRme}| = \frac{k\cdot \rme(k)}{\varphi(k)}. \]
\end{thm}
\noindent
\prf
Let us group the elements of \kRme\ into equivalence classes, according to the equivalence relation of Definition \ref{rn39}. By Theorem \ref{rn38}, we have that each equivalence class has $\varphi(k)$ elements, so the number of equivalence classes is
$\frac{|\kRme|}{\varphi(k)}$. Each representative of an equivalence class, has an orbit consisting of $k$ elements, so
we see that the above relation holds. \sqr

\begin{thm} \label{fs06}
If $m$ is weakly even, then $k\mapsto |\orbm{\kRme}|\in\MU$ for all $e\in\Em$.
\end{thm}
\noindent
\prf
Our theorem follows easily from Theorem \ref{fs02} and the above relation. \sqr

\newpage
\section{Idempotent Numbers as an Algebraic Structure}

\begin{dfn} \label{ia01}
For $e,e_1,e_2\im\Em$, let us define the following operators
\[ \be := (1-e)\ \mmod m \]
\[ e_1\co e_2 := (e_1 e_2+\beo \bet)\ \mmod m \]
\[ e_1\st e_2 := \overline{\beo\cdot\bet} \]
\[ e_1\sim e_2 := \overline{\beo\cdot e_2}. \]
\end{dfn}

\begin{thm} \label{ia02}
For $e\in\Em,\ a,b,c,d\in\Z,\ n\in\N$, the following identities hold
\[ (ae+b\be)(ce+d\be)\equiv (ac)e+(bd)\be\ (\mmod m) \]
\[ (ae+b\be)^n\equiv (a^n)e+(b^n)\be\ (\mmod m). \]
\end{thm}
\noindent
\prf
\[ (ae+b\be)(ce+d\be)\equiv (ac)e+(ad)e\be+(bc)\be e+(bd)\be\equiv (ac)e+(bd)\be\ (\mmod m). \]
The second identity follows from the first one. \sqr

\begin{thm} \label{ia03}
For $e,e_1,e_2\in\Em$, we have
\[ \be,\ e_1\co e_2,\ e_1\st e_2,\ e_1\sm e_2\in\Em. \]
\end{thm}
\noindent
\prf
Follows trivially from our previous theorem. \sqr

\begin{dfn} \label{ia04}
Let $\B_m$ denote the set
\[ \B_m:= \{\piai:\ \alpha_i>0\}. \]
For $A\subset\B_m$
\[ \bar{A}:= \B_m\setminus A. \]
For $e\in\Em$ define
\[ \B_m(e):= \{k\in\B_m:\ k\mid e\}. \]
\end{dfn}

\begin{thm} \label{ia05}
For $e,e_1,e_2\in\Em$ the following identities hold
\[ \B_m(\be)=\overline{\B_m(e)} \]
\[ \B_m(e_1\cdot e_2)=\B_m(e_1)\cup\B_m(e_2) \]
\[ \B_m(e_1\st e_2)=\B_m(e_1)\cap\B_m(e_2) \]
\[ \B_m(e_1\sm e_2)=\B_m(e_1)\setminus\B_m(e_2) \]
\[ \B_m(e_1\co e_2)=\B_m(e_1)\sd\B_m(e_2). \]
\end{thm}
\noindent
\prf
The first two identities are trivial.
\[ \B_m(e_1\st e_2) = \B_m(\overline{\beo\cdot\bet}) = \overline{\B_m(\beo\cdot\bet)} =
\overline{\overline{\B_m(e_1)}\cup\overline{\B_m(e_2)}} = \B_m(e_1)\cap\B_m(e_2) \]
\[ \B_m(e_1\sm e_2) = \B_m(e_1)\cap\overline{\B_m(e_2)} = \B_m(e_1)\setminus\B_m(e_2). \]
To prove the last identity, we first make a bit of calculation.
\[ (e_1\sm e_2)\cdot (e_2\sm e_1)\equiv \overline{\beo e_2}\cdot \overline{e_1\bet}\equiv (1-\beo e_2)(1-e_1\bet)\equiv \]
\[ \equiv 1-e_1\bet-\beo e_2\equiv e_1+\beo-e_1\bet-\beo e_2\equiv e_1\co e_2\ (\mmod m) \]
\[ \B_m(e_1\co e_2) = \B_m((e_1\sm e_2)\cdot (e_2\sm e_1)) = \B_m(e_1\sm e_2)\cup\B_m(e_2\sm e_1) = \]
\[ = (\B_m(e_1)\setminus\B_m(e_2))\cup (\B_m(e_2)\setminus\B_m(e_1)) = \B_m(e_1)\sd\B_m(e_2).\ \ \sqr \]

Considering the properties above, we see that an isomorphism may be defined between $\E_m$ and the class of subsets of any finite set, which has $\omega(m)$ elements.

\begin{thm} \label{ia06}
The structure $\langle \Em; \{1,^{-1},\co\} \rangle$ is an Abelian group, where each element is of order two.
\end{thm}
\noindent
\prf
First, we will show that for all $e_1,e_2\in\Em$, there exists one and only one $e_3\in\Em$ such that $e_1=e_2\co e_3$. Let us take any $e\in\Em$. Then
\[ e_2\co e\equiv \bet-(\bet-e_2)e\ (\mmod m) \]
and $(\bet-e_2,m)=1$ since $(\bet-e_2)^2\equiv 1\ (\mmod m)$. So by Theorem \ref{in07} we have
\[ |e_2\co\Em| = |\bet-(\bet-e_2)\Em| = |(\bet-e_2)\Em| = 2^{\omega\left(\frac{m}{(\bet-e_2,m)}\right)} = 2^{\omega(m)} = |\Em| \]
which proves both the existence and unicity of $e_3$. It is clear that $e\co 1=e$ and $e\co e=1$, so we have the existence of an inverse, and that each element is of order two. It is also obvious that $\co$ is commutative. In order to show that $\co$ is associative, take any $e_1,e_2,e_3\in\Em$. Then
\[ (e_1\co e_2)\co e_3 \equiv (e_1 e_2+\beo\bet)\co e_3\equiv (e_1 e_2+\beo\bet)e_3 + (\beo e_2+e_1\bet)\beh \equiv \]
\[ \equiv e_1 e_2 e_3 + \beo\bet e_3 + \beo e_2 \beh + e_1 \bet \beh \equiv e_1(e_2 e_3+\bet\beh) + \beo (\bet e_3 + e_2\beh) \equiv \]
\[ \equiv e_1\co (e_2 e_3 + \bet\beh) \equiv e_1\co(e_2\co e_3)\ (\mmod m).\ \ \sqr \]

\begin{thm} \label{ia07}
The $\st$ operator is commutative and associative. Multiplication is distributive with respect to $\st$, and $\st$ is distributive with respect to $\co$.
\end{thm}
\noindent
\prf
The commutativity of $\st$ is trivial. Now, take any $e_1,e_2,e_3\in\Em$.
\[ (e_1\st e_2)\st e_3 \equiv \overline{\overline{(e_1\st e_2)}\cdot\beh} \equiv
\overline{\overline{(\overline{\beo\cdot\bet})}\cdot\beh} \equiv \overline{\beo\cdot\bet\cdot\beh}\ (\mmod m) \]
which proves the associativity of $\st$.\\
Now, to prove the third property, we calculate
\[ e_1\cdot (e_2\st e_3)\equiv e_1\cdot\overline{\bet\cdot\beh}\equiv e_1(1-(1-e_2)(1-e_3))\equiv e_1e_3 + e_1e_2 - e_1e_2e_3 \equiv \]
\[ \equiv 1-(1-e_1e_2)(1-e_1e_3) \equiv \overline{\overline{e_1e_2}\cdot\overline{e_1e_3}}
\equiv (e_1\cdot e_2)\st(e_1\cdot e_3)\ (\mmod m). \]
The fourth property follows from
\[ e_1\st (e_2\co e_3)\equiv \overline{\beo\overline{(e_2e_3+\bet\beh)}}\equiv \overline{\beo(e_2\beh+\bet e_3)}\equiv \]
\[ \equiv \overline{\beo\beh-\beo\bet\beh+\beo\bet-\beo\beh\bet}\equiv \overline{\overline{\beo\bet}\cdot \beo\beh +
 \beo\bet\cdot\overline{\beo\beh}}\equiv \]
\[ \equiv \overline{\beo\bet}\cdot\overline{\beo\beh}+\beo\bet\cdot\beo\beh\equiv (e_1\st e_2)\co (e_1\st e_3)\ (\mmod m).\ \ \sqr \]

\begin{thm} \label{ia08}
The structure $\langle \Em; \{\co,\st\} \rangle$ is a commutative ring, with $\co$ being ``addition'' and $\st$ being ``multiplication''.
\end{thm}
\noindent
\prf
Follows from our previous two theorems. \sqr

We see that because of the isomorphism that exists between $\Em$ and the subsets of a finite set, the above theorem states the
well-known fact from Set Theory, that the subsets of a set form a commutative ring with respect to the operators $\cap$ and $\sd$.

\begin{thm} \label{ia09}
For $e,e_1,e_2\in\Em$, we have
\[ e\cdot\be\equiv m,\ e+\be\equiv 1\ (\mmod m) \]
\[ e\co 1=e,\ e\co\be=m,\ e\co m=\be \]
\[ \overline{e_1\co e_2}=\beo\co e_2=e_1\co\bet \]
\[ e_1\co e_2\equiv (e_1+\bet)(\beo +e_2)\equiv (e_1-\bet)^2\equiv (\beo-e_2)^2\ (\mmod m). \]
\end{thm}
\noindent
\prf
The fourth line of identities seems a bit nontrivial, so we shall prove it in part below.
\[ (e_1+\bet)(\beo+e_2)\equiv e_1\beo + e_1 e_2 + \bet\beo + \bet e_2\equiv e_1 e_2 + \bet\beo\ (\mmod m) \]
\[ (e_1-\bet)^2 \equiv e_1-2e_1\bet +\bet \equiv e_1(1-\bet)+\bet(1-e_1)\equiv e_1 e_2 + \bet\beo\ (\mmod m).\ \ \sqr \]

Note that the first and second lines of identities show that $\co$ behaves somewhat like multiplication. In our upcoming theorems, we will prove properties of $\st$ which show that it may behave in a sense both like multiplication and addition.

\begin{thm} \label{ia10}
For $e,e_1,\dots,e_n\in\Em$, we have
\[ e\st e = e,\ e\st 1 = 1 \]
\[ (e_1\st e_2)-(\beo\st\bet)\equiv e_1\cdot e_2-\beo\cdot\bet\ (\mmod m) \]
\[ ((e_1\st e_2)-(\beo\st\bet))^2\equiv e_1\co e_2\ (\mmod m) \]
\[ \overline{\beo\st\bet} \equiv e_1\cdot e_2\ (\mmod m) \]
\[ \bigotimes_{i=1}^n e_i = \overline{\prod_{i=1}^n \bar{e}_i} \]
\[ e\st\be=1,\ e\st 0=e \]
\[ (e_1\cdot e_2)\st (\beo\cdot\bet)\equiv e_1\cdot e_2+\beo\cdot\bet\equiv e_1\co e_2\ (\mmod m). \]
\end{thm}
\noindent
\prf
The first four lines of properties, are quite trivial. The fifth property may be proven via induction, using the associativity of $\st$. The sixth and seventh lines are quite trivial calculations as well. \sqr

Our next theorem shows a peculiar property of $\st$, in which it behaves both like addition and multiplication. In fact, the second property sheds light on the double nature of this operator.

\begin{thm} \label{ia11}
For $e,e_1,e_2\in\Em$, we have
\[ (e_1\co e)\st (e_2\co e) \equiv (e_1\st e_2)e + (\beo\st\bet)\be\ (\mmod m) \]
\[ e_1\st e_2\equiv e_1+e_2 -e_1\cdot e_2\ (\mmod m). \]
\end{thm}
\noindent
\prf
\[ (e_1\co e)\st (e_2\co e)\equiv 1-\overline{e_1\co e}\cdot\overline{e_2\co e}\equiv \]
\[ \equiv 1-(\beo e+e_1\be)(\bet e+e_2\be)\equiv 1-(\beo\bet e+e_1 e_2 \be)\equiv \]
\[ \equiv e+\be - (\beo\bet e+e_1 e_2 \be)\equiv (e_1\st e_2)e + (\beo\st\bet)\be\ (\mmod m). \]
The second property is just simple calculation. \sqr

\newpage
\section{Second-Degree Polynomials}

\begin{dfn} \label{sd01}
For $k\in\Z$, let $\Smk$ denote the set of solutions of the equation
\[ x^2\equiv k x\ (\mmod m) \]
among the elements of $\Zm$.
\end{dfn}

\begin{thm} \label{sd02}
Let $k\in\Z$ be such that $(k,m)=1$. Then
\[ \Smk = k\Em\ \mmod m. \]
\end{thm}
\noindent
\prf
Let $i\in\N$ be such that $\ali$ is positive. Then from
\[ x_0^2\equiv kx_0\ (\mmod \piai) \]
it follows that
\[ x_0\equiv 0\ \mor\ x_0\equiv k\ (\mmod \piai). \]
Let $e\in\Em$ be such that $\mu_m(x_0)=\mu_m(e)$. Then for all $i\in\N$ we have $x_0\equiv ke\ (\mmod \piai)$, so $x_0\equiv ke\ (\mmod m)$. So we may conclude that $\Smk\subset k\Em\ \mmod m$.\\
Now, we see that $k\Em\ \mmod m\subset\Smk$ as well, since for any $e\in\Em$, we have
\[ (ke)^2\equiv k(ke)\ (\mmod m).\ \ \sqr \]

\begin{thm} \label{sd03}
Let $a,b\in\Z$ be such that $(b-a,m)=1$. Then for all solutions $r\in\Z_m$ of the equation
\[ (x-a)(x-b)\equiv 0\ (\mmod m) \]
there exists a unique $e\in\Em$, such that
\[ r\equiv ae+b\be\ (\mmod m). \]
\end{thm}
\noindent
\prf
Our equation may be rearranged as
\[ (x-a)^2\equiv (b-a)(x-a)\ (\mmod m). \]
So since $(b-a,m)=1$, by our previous theorem we have that for all solutions $r\in\Z$, there exists a unique $e\in\Em$, such that
\[ r-a\equiv (b-a)e\ (\mmod m) \]
which may be rearranged as
\[ r\equiv ae+b\be\ (\mmod m).\ \ \sqr \]

\begin{thm} \label{sd04}
Let $m$ and $a\in\Z$ be such that $(2a,m)=1$. Then for all solutions $r_1,r_2\in\Z_m$ of the equation
\[ x^2\equiv a\ (\mmod m) \]
there exists a unique $e\in\Em$, such that
\[ r_1\equiv r_2(e-\be)\ (\mmod m). \]
\end{thm}
\noindent
\prf
With the notation above, we have that our equation is equivalent to the equation
\[ (x-r_2)(x-(-r_2))\equiv 0\ (\mmod m). \]
Now, since $r_2^2\equiv a\ (\mmod m)$, we have $(r_2,m)=(a,m)=1$, from which we have $(r_2-(-r_2),m)=1$ since $2\nmid m$, so by our previous theorem, we have that there exists a unique $e\in\Em$, such that
\[ r_1\equiv r_2 e+(-r_2)\be\equiv r_2(e-\be)\ (\mmod m).\ \ \sqr \]

\begin{thm} \label{sd05}
Let $m$ be an odd number, or four times an odd number. Then for all solutions $r\in\Zm$ of the equation
\[ x^2\equiv 1\ (\mmod m) \]
there exists a unique $e\in\Em$, such that
\[ r\equiv e-\be\ (\mmod m). \]
\end{thm}
\noindent
\prf
The case when $m$ is odd, follows from our previous theorem. Now, if $m$ is four times an odd number, then it is easy to see that
\[ \omega\left(\frac{m}{(2,m)}\right) = \om. \]
Our equation is equivalent to the equation
\[ (x+1)^2\equiv 2(x+1)\ (\mmod m) \]
so we see that all elements of $2\Em-1$ satisfy this equation, and by Theorem \ref{in07} we also have that
\[ |2\Em-1\ \mmod m|=2^{\omega\left(\frac{m}{(2,m)}\right)}=2^{\om} \]
which is the number of solutions of our equation if $m$ is four times an odd number, so we have that, for all $r\in\Zm$ satisfying the equation, there exists a unique $e\in\Em$, such that
\[ r\equiv 2e-1\equiv e-\be\ (\mmod m).\ \ \sqr \]

\begin{thm} \label{sd15}
Let $m$ be an odd number. Then for all $e\in\Em$
\[ \SmR(2,e)\subset \{e(e_0-\be_0)\ \mmod m: e_0\in\Em\}. \]
Furthermore, if $e\neq m$, then for all $e_0\in\Em$ we have $e(e_0-\be_0)\not\equiv e(\be_0-e_0)\ (\mmod m)$. Moreover, the following properties are valid
\[ |\SmR(2,e)| = 2^{\omega(\mu_m(e))} \]
\[ \prod \SmR(2,e)\equiv (-1)^{2^{\omega(\mu_m(e))-1}}\cdot e\ (\mmod m). \]
\end{thm}
\noindent
\prf
Take any $i\in\N$ such that $p_i\neq 2$, and $a\in\SmR(2,e)$. There are two possible cases. If $a^2\equiv 0\ (\mmod \piai)$, then $a\equiv 0\ (\mmod \piai)$, since $a\in\R_{\piai}$. If $a^2\equiv 1\ (\mmod \piai)$, then since $p_i\neq 2$ it follows that
$a\equiv \pm 1\ (\mmod \piai)$. From these two cases, we have that $a\equiv e(e_0-\be_0)\ (\mmod m)$, for some $e_0\in\Em$.\\
Let us suppose indirectly, that there exists some $e_0\in\Em$, such that $e(e_0-\be_0)\equiv e(\be_0-e_0)\ (\mmod m)$. Then $2(e_0-\be_0)e\equiv 0\ (\mmod m)$, from which we have $e\equiv 0\ (\mmod m)$, since $(2(e_0-\be_0),m)=1$, because
\[ (e_0-\be_0)^2\equiv e_0+\be_0\equiv 1\ (\mmod m). \]
So $e=m$, which of course is a contradiction.\\
To prove the third property, observe that for all $a\in\SmR(2,e)$
\[ a^2\equiv e\ (\mmod m)\ \lrar\ a^2\equiv 1\ (\mmod \mu_m(e))\ \mand\ a^2\equiv 0\ \left(\mmod \frac{m}{\mu_m(e)}\right). \]
So for all $\piai\in\overline{\Bm(e)}$, we have $a\equiv\pm 1\ (\mmod \piai)$. Meanwhile $a\equiv 0\ (\mmod \frac{m}{\mu_m(e)})$, since $a\in\Rme$. So by the Chinese Remainder Theorem, and since $\omega(\mu_m(e))=|\overline{\Bm(e)}|$, we have the formula.\\
By the first property, with the notation $n:=2^{\omega(\mu_m(e))}$, we have
\[ \prod \SmR(2,e) \equiv e(e_1-\be_1)(e_2-\be_2)\dots (e_n-\be_n)\equiv \]
\[ \equiv e\cdot\left[(e_1-\be_1)\dots (e_{\frac{n}{2}}-\be_{\frac{n}{2}})\right]\cdot \left[(\be_1-e_1)\dots (\be_{\frac{n}{2}}-e_{\frac{n}{2}}) \right] \equiv \]
\[ \equiv e\cdot (-1)^{\frac{n}{2}}\ (\mmod m).\ \ \sqr \]

\begin{dfn} \label{sd06}
For $k\in\Z,\ r\in\Smk,\ e\in\Em$ define
\[ \br:= (k-r)\ \mmod m \]
\[ r\co e:= (re+\br\be)\ \mmod m \]
\[ r\st e:= (k-\br\be)\ \mmod m. \]
\end{dfn}

Note that we continue to use the same notations as in the previous section. In order to distinguish between these operators that have been denoted the same way, even though they are different, always refer to the set from which the operands have been taken.

\begin{thm} \label{sd07}
For $k\in\Z,\ r\in\Smk,\ e\in\Em$, we have
\[ \br,\ r\co e,\ r\st e\in\Smk. \]
\end{thm}
\noindent
\prf
\[ (k-r)^2\equiv k^2-2kr+r^2\equiv k^2-2kr+kr\equiv k(k-r)\ (\mmod m) \]
\[ (r\co e)^2\equiv r^2 e+\br^2\be\equiv (kr)e+(k\br)\be\equiv k(r\co e)\ (\mmod m) \]
\[ (r\st e)^2\equiv k^2-2k\br\be+\br^2\be\equiv k(k-\br\be)\equiv k(r\st e)\ (\mmod m).\ \ \sqr \]

\begin{thm} \label{sd08}
Take any $e\in\Em$ and $k\in\Z$. Then
\[ \Smk = \Smk\co e. \]
\end{thm}
\noindent
\prf
In order to show the equality of the two sets, it is enough for us to prove that for any $r_1,r_2\in\Smk,\ r_1\neq r_2$, we have $r_1\co e\neq r_2\co e$. For let us suppose indirectly, that there exist some $r_1,r_2\in\Smk,\ r_1\neq r_2$, such that
$r_1\co e = r_2\co e$. Then
\[ r_1-r_2\equiv (e-\be)^2 (r_1-r_2)\equiv (e-\be)((r_1-r_2)e+(r_2-r_1)\be)\equiv \]
\[ \equiv (e-\be)((r_1-r_2)e+(\br_1-\br_2)\be)\equiv (e-\be)(r_1\co e - r_2\co e)\equiv 0\ (\mmod m). \]
So we arrive at a contradiction. \sqr

\begin{thm} \label{sd13}
For $k\in\Z,\ e,e_1,e_2\in\Em,\ r\in\Smk$ the following properties hold
\[ \overline{r\co e} = r\co\be = \br\co e \]
\[ (r\co e_1)\co e_2 = r\co (e_1\co e_2). \]
\end{thm}
\noindent
\prf
\[ \overline{r\co e}\equiv k-(r\co e)\equiv ke+k\be-(re+\br\be)\equiv \br e+r\be\equiv r\co\be\equiv \br\co e\ (\mmod m) \]
\[ (r\co e_1)\co e_2\equiv (re_1+\br\be_1)\co e_2\equiv (re_1+\br\be_1)e_2+(\br e_1+r\be_1)\be_2\equiv \]
\[ \equiv re_1 e_2+\br\be_1 e_2+\br e_1\be_2+r\beo\bet\equiv r(e_1e_2+\beo\bet)+\br(\beo e_2+e_1\bet)\equiv \]
\[ \equiv r\co(e_1\co e_2)\ (\mmod m).\ \ \sqr \]

Note that the second property is somewhat like associativity.

\begin{thm} \label{sd14}
For $k\in\Z,\ (k,m)=1$ and $r\in\Smk$, the following equivalence holds
\[ r\co e_1 = r\co e_2\ \lrar\ e_1=e_2\ \ (e_1,e_2\in\Em). \]
\end{thm}
\noindent
\prf
The $\lar$ part of the equivalence is trivial. To prove the $\rar$ part, first we see that
\[ (r-\br)^2\equiv r^2-2r\br+\br^2\equiv kr+k\br\equiv k^2\ (\mmod m). \]
Now, let us suppose that $r\co e_1 = r\co e_2$. This means that
\[ re_1+\br\be_1\equiv re_2+\br\be_2\ \rar\ \br+(r-\br)e_1\equiv \br+(r-\br)e_2\ (\mmod m). \]
First subtracting $\br$, then squaring both sides, we get
\[ k^2 e_1\equiv k^2 e_2\ (\mmod m) \]
which implies that $e_1=e_2$, since $(k,m)=1$. \sqr

\begin{dfn} \label{sd09}
For $e\in\Em,\ r_1,r_2\in\Sy_{m,e}$ define
\[ r_1\co r_2 := (r_1r_2+\br_1\br_2)\ \mmod m \]
\[ r_1\st r_2 := \overline{\br_1\cdot\br_2}. \]
\end{dfn}

Same can be said for these operators, as for those of Definition \ref{sd06}.

\begin{thm} \label{sd10}
For $e\in\Em,\ r_1,r_2\in\Sy_{m,e}$, we have
\[ r_1\co r_2,\ r_1\st r_2\in\Sy_{m,e}. \]
\end{thm}
\noindent
\prf
\[ (r_1\co r_2)^2\equiv r_1^2 r_2^2 + \br_1^2 \br_2^2\equiv e(r_1\co r_2)\ (\mmod m) \]
\[ (\overline{\br_1\cdot\br_2})^2\equiv e-2e\br_1\br_2+ \br_1^2 \br_2^2\equiv e(1-\br_1\br_2)\equiv e(r_1\st r_2)\ (\mmod m).\ \ \sqr \]

\begin{thm} \label{sd11}
For $e\in\Em,\ a,b\in\Rme,\ c,d\in\Z,\ n\in\N,\ r\in\Sy_{m,e}$, we have
\[ (ar+b\br)(cr+d\br)\equiv (ac)r+(bd)\br\ (\mmod m) \]
\[ (ar+b\br)^n\equiv a^n r+b^n\br\ (\mmod m). \]
\end{thm}
\noindent
\prf
\[ r\br\equiv 0\ \rar\ (ar+b\br)(cr+d\br)\equiv (ac)r^2 +(bd)\br^2\equiv (ea)cr+(eb)d\br\ (\mmod m). \]
The second property follows from the first one via induction. \sqr

\begin{thm} \label{sd12}
For all $e\in\Em,\ k\in\Rme$ we have
\[ \Smk\cap\Rme = \{k\}. \]
\end{thm}
\noindent
\prf
For any $r\in\Smk\cap\Rme$ we have
\[ r\in_{\mu_m(e)} \Sy_{\mu_m(e),k}\cap \R_{\mu_m(e)}^1 = (k\E_{\mu_m(e)}\ \mmod \mu_m(e))\cap \R_{\mu_m(e)}^1 = \]
\[ = \{k\ \mmod \mu_m(e)\}\ \rar\ r\equiv k\ (\mmod \mu_m(e)). \]
So, since
\[ r\equiv 0\equiv k\ \left(\mmod \frac{m}{\mu_m(e)}\right) \]
and $\mu_m(e)=\mu_m(k)=\mu_m(r)$, we have that $r=k$. \sqr

\newpage
\phantomsection
\section*{Acknowledgements}
\addcontentsline{toc}{section}{Acknowledgements}

I would like to express my vast gratitude to Professor Mih\'aly Szalay, for the enormous amount of time and effort it took to review my paper. Besides his various helpful comments and corrections, I would like to say special thanks to him for his help on the proof of Theorem \ref{rn15}, which was a lemma of utmost importance.


\newpage
\bibliographystyle{abbrv}
\bibliography{mybib}
\addcontentsline{toc}{section}{\textbf{References}}

\end{document}